\documentclass[11pt,fleqn]{article}
\pdfoutput=1
\usepackage{amsmath}
\usepackage{amssymb}
\usepackage{makeidx}
\usepackage{pdfpages}
\usepackage{pgf,tikz}
\usepackage{mathrsfs}
\usetikzlibrary{arrows}
\usepackage{amssymb,amsmath,graphicx}
\usepackage{epstopdf}
\usepackage{multirow}
\usepackage{caption}
\usepackage{blkarray}

\begin{document}
\sloppy

\newtheorem{axiom}{Axiom}[section]
\newtheorem{example}[axiom]{Example}
\newtheorem{proposition}[axiom]{Proposition}

\newcommand{\eps}{\varepsilon}
\newcommand{\real}{\mathbb{R}}
\newcommand{\nat}{{\mathbb{N}}}
\newcommand{\seq}[1]{\langle #1\rangle}

\newcommand{\crcl}{\mbox{$\mathcal{R}$}}
\newcommand{\bound}{\mbox{$\mathcal{L(R)}$} }
\newcommand{\ccc}{{\cal C}}
\newcommand{\VL}{\mbox{$\mathrm{VL}$}} 
\newcommand{\VR}{\mbox{$\mathrm{VR}$}}


\title{2VRP: a benchmark problem  for small but rich VRPs.}

\author{
  Vladimir Deineko
  \thanks{{\tt V.Deineko@warwick.ac.uk}.
  Warwick Business School, Coventry, CV4 7AL, United Kingdom}
  \thanks{Corresponding author; tel.\ +44 02476524501}
  \and Bettina Klinz
  \thanks{{\tt klinz@opt.math.tu-graz.ac.at}.
Institut f\"ur Diskrete Mathematik, TU Graz, Steyrergasse 30, A-8010 Graz, Austria}
}
\date{}
\maketitle

\maketitle



\begin{abstract}
We consider a 2-vehicle routing problem (2VRP) which can be viewed as a building block for the variety of vehicle routing problems (VRP). As a simplified version of the 2VRP, we consider a 2-period balanced travelling salesman problem (2TSP) and describe a polynomially solvable case of this ${\cal NP}$-hard problem. For the 2VRP with  general settings, we suggest a framework based on the Held and Karp dynamic programming algorithm. Our algorithms based on this framework show an exceptionally good performance on the published test data. Our approach can be easily extended to a variety of constraints/attributes in the VRP, hence the wording ``small but rich" in the title of our paper. We also introduce a new methodological approach: we use easy solvable special cases for generating test instances and then use these instances in computational experiments.

\end{abstract}

\section{Introduction and related work}
\label{sec:introduction}
\nopagebreak

In the vehicle routing problem (VRP) a set of customers with  certain  requests are to be visited by vehicles. The vehicles are to be chosen from a fleet of heterogeneous vehicles, with various fixed and variable costs of usage. The objective is to find a minimum cost schedule for customer visits to deliver required services in the specified time and manner. 

It would not be an exaggeration to say that thousands of research papers devoted to the VRPs are published every year. As rightly mentioned by Michael Drexl \cite{Drexl}, this ``is certainly due to the intellectual challenge VRPs pose as well as to their relevance in logistics and transport''.

The VRP can be viewed as a combination of the two well known combinatorial optimization problems - the travelling salesman problem (TSP) (see e.g. \cite{ABCC}) and the bin packing/knapsack problem (see e.g. \cite{KPP}). It is not surprising that the combination of these two ${\cal NP}$-hard problems creates new computational challenges for researchers and practitioners. For instance, the sizes of the TSP instances which are tractable by recently developed computational algorithms \cite{ABCC} are much bigger than the sizes of easy tractable VRP instances.

\smallskip
{\bf Two vehicle routing problem as an atomic model for the VRP.} In many practical applications of the VRP, the number of customers visited by a single vehicle is not very large. We mention here just a few examples: teachers visiting special needs pupils, nurses attending patients at home, food delivery in rural areas (long distances - hence few customers to visit), bulk deliveries of industrial goods (large items - hence few customers). On the other hand the myriads of real life  constraints make even small size problems computationally challenging.
In recent publications, the VRPs with many practice related constraints are referred to  as \emph {rich} VRPs (see surveys \cite{CAGR,Drexl,LKH}), or \emph{multi-attribute} VRPs \cite{VCGP1,VCGP2}.
In our study we focus on the simplest possible version of the rich VRP - the VRP with many constraints/attributes but with only two vehicles, which we call the 2VRP. Similar to the VRP with many vehicles, the 2VRP encompasses hard optimisation problems mentioned above, the TSP and the knapsack problem. In some sense, the 2VRP can be viewed as an atomic model for the VRP.

Further on, keeping in mind the practical applications mentioned above, we studied small (with not many customers) but rich VRPs (hence the title of the paper). We hope that better understanding of the 2VRP would allow researchers to generate new ideas, and hence to construct better algorithms for the rich VRPs with many vehicles. 

\bigskip
{\bf Two-period travelling salesman problem.}
We start our study of the 2VRP with an investigation of 
a simplified version of this problem, which is also known as the 2-period TSP (2TSP). In the classic TSP, a salesman is given an $n\times n$ distance matrix among $n$ locations: location $1$ is where the salesman lives, and $n-1$ locations are for the customers to be visited. The objective of the salesman is to find the shortest route to visit all customers and return home.
In the 2TSP version, a set of $m$ customers is to be visited in each of the two periods,  and $ k $ customers are to be visited once in either of the two periods. The objective is to achieve the minimal total distance travelled. 

Butler, Williams, and Yarrow \cite{ButW97} considered the 2TSP as a model for a milk collection problem in Ireland. They applied what they called a ``man-machine method'', combining an integer programming technique with a human being intervention for identifying violated constraints. 

\begin{figure}
\unitlength=1cm
\begin{center}
{\begin{picture}(6.5,6)
\includegraphics[scale=1.1]{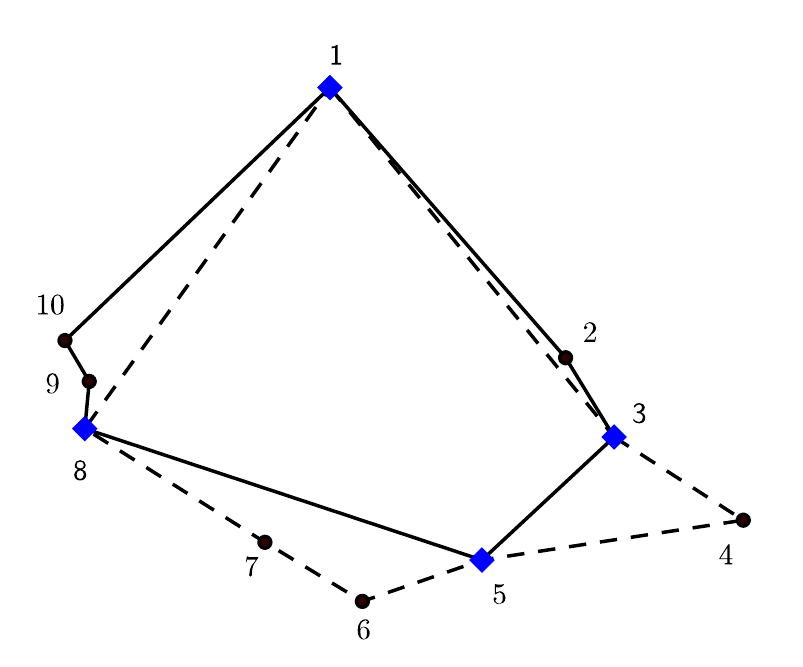}
\end{picture} }
\end{center}
\caption{Illustration to the definition of the balanced 2TSP; 
customers $1$, $3$, $5$, and $8$ are visited in two periods; the same number of customers is visited in each of the two periods.}
\label{fig:2pTSP}
\end{figure}

In the \emph{balanced} 2TSP an additional constraint requests that the number of customers visited in each period differs by no more than 1.
Bassetto and Mason \cite{BassettoM11} considered a balanced 2TSP in the Euclidean plane. In this variant of the problem the customer locations are points in the Euclidean plane, and the costs of travelling between the customers are standard Euclidean distances. 
Fig.~\ref{fig:2pTSP} illustrates an instance of the Euclidean balanced 2TSP with $10$ customers. Four out of ten customers in this example are visited in $2$ periods (two sub-tours); the depot here is  a customer which is visited twice (see Fig.~\ref{fig:inst55}  for another instance of the 2TSP in the Euclidean plane).

The 2TSP can be modelled as a special case of the 2VRP  as follows. One of $m$ customers that have to be visited in two periods is chosen as a depot. For each of the other $m-1$ customers an identical copy (i.e. a new customer in the same location) is created. The identical customers are allocated then to two different vehicles. The allocation of these customers to the vehicles is fixed: each vehicle has to serve $m-1$ fixed customers, while $k$ customers to be served by only one of the vehicles. 

So, the balanced 2TSP can be viewed as the 2VRP with the following specific attributes:
\begin{itemize}
\item There are $(k+2\times (m-1))$ customers with unit demand for  deliveries and one depot. For simplification of notations we assume that $k$ is even.
\item There are two homogeneous vehicles with the capacity of 
$m-1+ k/2 $ which are located in the same depot.
\item Each of the vehicles has $(m-1)$ \emph{fixed} customers which are 
pre-allocated for visits by the vehicle;
\item The objective is to find the routes for the vehicles with the total minimal distance travelled, subject to the constraint that each vehicle visits $(m-1+ k/2)$ customers.
\end{itemize}

As the first step in our investigation, we describe a new polynomially solvable case of the balanced 2TSP. Our results are related to the well known solvable cases the standard TSP.

\bigskip

{\bf Polynomially solvable cases of the TSP.}
A feasible solution to the TSP can be represented as a  \emph{tour} $\tau=\seq{1=\tau_1,\tau_2,\ldots,\tau_n,1=\tau_{n+1}}$: a salesman starts from location $1$, visits then customer $\tau_2$, and so on, and eventually returns to the initial location.  We will refer to the items in tours as \emph{customers}, or \emph{cities}, or \emph{nodes}. 
Given a distance matrix $(c_{ij})$, the length of tour $\tau$ is calculated as $c(\tau)=\sum_{i=1,\ldots,n}c_{\tau_i,\tau_{i+1}}$.

In the context of the VRP, algorithms for the TSP can be useful as procedures for improvements of routes for the vehicles.
Since the TSP is an  ${\cal NP}$-hard problem, investigation of polynomially solvable special cases of the TSP is the well established branch of research (see surveys \cite{BSurv,DKTW,GLS,Kabadi} for further references).
Probably one of the best-known special cases of the TSP is the case with the so-called \emph{Kalmanson} distance matrix.

A \emph{symmetric} $n\times n$ matrix $C$ is called a \emph{Kalmanson} matrix if it fulfils the 
\emph{Kalmanson conditions}:
\begin{eqnarray}
c_{ij}+c_{\ell m}\le c_{i\ell}+c_{jm}, && 
 \label{kalm1.c} \\
c_{im}+c_{j\ell}\le c_{i\ell}+c_{jm}, &&
\mbox{~for all~} 1\le i<j<\ell<m\le n. \label{kalm2.c}
\end{eqnarray}
Kalmanson~\cite{Kalmanson} noticed that if $n$ points in the Euclidean plane are located on the boundary of their convex hull and numbered along the convex hull, then the distance matrix for these points satisfies conditions (\ref{kalm1.c})-(\ref{kalm2.c}). Notice, that if a distance matrix satisfies conditions (\ref{kalm1.c})-(\ref{kalm2.c}), it does not necessarily mean that the points are on the boundary of a convex hull. For instance, the points in the example in Fig.~\ref{fig:2pTSP}, which we have adopted from \cite{DRVW}, are not on the boundary of the convex hull, however the distance matrix for these points is indeed the Kalmanson distance matrix (see Fig. 1 in \cite{DRVW} for the coordinates of the points in the Euclidean plane). 

The TSP with a Kalmanson distance matrix is solved to optimality 
by the tour $\pi=\seq{1,2,3,\ldots,n-2,n-1,n,1}$. It can be easily shown that cyclic renumbering of rows and columns does not destroy property (\ref{kalm1.c})-(\ref{kalm2.c}), so any permutation obtained as a cyclic shift of the identity permutation is also an optimal solution to the TSP.

Kalmanson matrices posses nice structural properties that allow one in polynomial time to find not only solutions to the TSP but also solutions to other well known combinatorial optimisation problems (see \cite{CFSS,DRW,DW,KW,PSW,Shmoys}). We single out here the so-called master tour problem (\cite{DRW}), which is useful  for our characterisation of a special solvable case of the balanced 2TSP (see Section~\ref{sec:solvable2TSP}). 

An  \emph{optimal} tour for the TSP is called the \emph{master tour}, if after deleting any subset of points from the tour, the so obtained sub-tour is still optimal for the remaining set of points. Given a distance matrix, \emph{the master tour problem} is to find out whether it is possible to construct the master tour. Surprisingly, the master tour exists if and only if the underlying distance matrix is a Kalmanson matrix (Theorem 5.1 in \cite{DRW}). In what follows we make use of two master tours, tour $\pi$ and the inverse of it, tour 
$\pi^-=\seq{1,n,n-1,\ldots,2,1}$.

\bigskip
A generalization of Kalmanson matrices is the class of Demidenko matrices.

A symmetric matrix $C = (c_{ij})$ is called a \emph{Demidenko matrix} if
\begin{eqnarray}\label{eq:demi}
c_{ij} + c_{\ell m} \leq c_{i\ell} + c_{jm},
&\mbox{\qquad for all $1\le i < j < \ell < m\le n$.}
\end{eqnarray} 

An optimal tour for the TSP with a Demidenko distance matrix can be found 
among so-called \emph{pyramidal} tours (see ~\cite{Demi} and \cite{GLS}). 
A tour $\tau=\seq{1,\tau_2,\ldots,\tau_m,n,\tau_{m+2},\ldots,\tau_{n},1}$ is called pyramidal, 
if $1<\tau_2<\ldots<\tau_m<n$, and $n>\tau_{m+2}>\ldots >\tau_{n}>1$.
We will also call a sub-tour obtained from $\tau$ by deleting any subset of items a pyramidal sub-tour. 

An optimal pyramidal tour can be found among $2^{n-2}$ tours in $O(n^2)$ time.  For the case of generality we consider here the case with an \emph{asymmetric} distance matrix. Length $E(1,1)$ of the optimal pyramidal tour can be found 
as $$E(1,1)=\min\{c_{12}+E(2,1),c_{21}+E(1,2)\},$$ where values $E(2,1)$ and $E(1,2)$ are calculated by using the following dynamic programming recursions:
\begin{gather*}
E(i,j)=\min
  \bigl\{ E(j+1,j)+c_{i,j+1},E(i,j+1)+c_{j+1,j} \bigr\},\\
E(j,i)=\min
  \bigl\{E(j+1,i)+c_{j,j+1},E(j,j+1)+c_{j+1,i} \bigr\},\\
\qquad i=1,2,\ldots,n-1; j=i+1,\ldots,n.\\
E(i,n)=c_{in};\ \  E(n,i)=c_{ni},\ \ 
\qquad i=1,\ldots,n.
\end{gather*}
Indeed, let $E(i,j)$, $i<j$, be the length of an optimal pyramidal sub-tour from $i$ to $j$ through the set of indices $\{j+1,j+2,\ldots,n\}$. Similarly, let $E(j,i)$, $i<j$, be the length of the shortest pyramidal sub-tour from $j$ to $i$.  Since we consider only pyramidal sub-tours, $j+1$ is placed either next to $i$ or next to $j$. In the optimal tour we chose the best option, i.e. the sub-tour with the minimal total length. This explains the recursions. 


In our computational experiments we use a TSP heuristic which uses the recursions above. Below we give a brief description of this heuristic.

\bigskip
{\bf A TSP heuristic that uses pyramidal tours.}
The algorithm for finding an optimal pyramidal tour searches an exponential neighbourhood in polynomial time. 
Renumbering customer locations (i.e.\ renumbering rows and columns in the distance matrix)  changes the neighbourhood searched by the algorithm.  Given a permutation, Carlier \& Villon's \cite{CV87} considered a large neighbourhood that contains sets of pyramidal tours constructed from all cyclic shifts of the initial permutation. The optimal permutation in this large neighbourhood is used as a new solution to the TSP and as a permutation for renumbering customers. The search is repeated until no better solution is found. 

We illustrate this idea on a small instance. Assume that $5\times 5$ distance matrix $C=(c_{ij})$ and an initial permutation  $\sigma_0=(1,2,3,4,5)$ are given. For permutation $\sigma_0$ and each of its cyclic shifts, i.e. for $\sigma_1=(2,3,4,5,1)$, $\sigma_2=(3,4,5,1,2)$, $\sigma_3=(4,5,1,2,3)$ and $\sigma_4=(5,1,2,3,4)$, an optimal pyramidal tour for the TSP with the permuted matrix $(c_{\sigma_k(i)\sigma_k(j)})$, $k=0,1,\ldots,4$, is found. $BELPERM(\sigma_0)$ returns tour $\tau^*$ which is the shortest tour among  $5$ tours found. The procedure is restarted then as $BELPERM(\tau^*)$.
 
Heuristic $BELPERM$ performed very well in computational experiments reported in \cite{CV87}. It also has nice theoretical properties: 
\begin{itemize}
\item  The heuristic returns a solution that cannot be approved by the well known $2$-optimal TSP heuristic \cite{CV87};
\item The large neighbourhood, which is searched on each iteration of the local search described above, contains at least $75\%$ of permutations searched by the  well-known $3$-optimal TSP heuristic (see Theorem 5 in \cite{Neighborhoods}). 
\end{itemize}

In our computational experiments, we use $BELPERM$ as a heuristic for improving individual routes in the 2VRP solutions. 

\bigskip
{\bf Computational experiments and benchmark problems.}
A possible approach to empirical evaluation of algorithms for hard optimisation problems is an extensive computational testing of the algorithms on published benchmark problems. 
For the 2TSP problem, conducting wide range computational experiments is problematic, since there are no many publications on this topic.  To the best of our knowledge, the only paper which reports computational experiments for the 2TSP 
is Bassetto \& Mason \cite{BassettoM11} paper. 
Therefore we decided to generate additional benchmark problems by using polynomially solvable cases of the 2TSP.

While there are some known polynomially solvable cases of the VRP (see  \cite{Archetti1, Archetti2, Hassin, Labbe, YuLiu}), we are not aware of any polynomially solvable case of the balanced 2VRP (2TSP). Recall that ``2" in the abbreviation 2VRP means the number of vehicles, not the capacity of the vehicle as for instance in the abbreviation 2SDVRP for the Split Delivery Vehicle Routing Problem \cite{Archetti1} (see also \cite{Hassin}, where 2 means the number of customers visited).

In this paper we describe a polynomially solvable case of the balanced 2TSP with Kalmanson matrices. We also describe an algorithm for generating random Kalmanson matrices, and use it then in our computational experiments. Notice, the distance matrices for the points located on a line, or on a cycle, or on a tree (the cases considered in \cite{Archetti1, Archetti2, Hassin, Labbe, YuLiu}) are all Kalmanson matrices. So our approach can also be used for testing heuristics developed, in particular, for the VRPs with the ``tree-like" networks (see \cite{Basnet, Chandran, Labbe}). 

\bigskip
{\bf Our contribution and structure of the paper.} 
In Section~\ref{sec:solvable2TSP}, we focus on a special case of the 2VRP which is also known as the balance 2-period TSP (2TSP). We describe a new polynomially solvable case of this problem and explain dynamic programming recursions that can be used for finding optimal solutions in this special case. We use these recursions in a heuristic for finding approximate solutions for the 2TSP in general case. Computational experiments on the benchmark problems from \cite{BassettoM11} show that this heuristic is competitive with the previously published algorithms. We also describe an algorithm for generating instances of the 2TSP with known optimal solutions. 

Our main contribution is a framework for rich 2VRPs, which we present in Section~\ref{sec:framework}. 
Computational tests on published 2TSP benchmark problems have shown exceptionally good performance  of our framework. To further investigate the efficiency of our approach we conducted extensive computational experiments on the set of instances with known optimal solutions, which were generated as described in Section~\ref{sec:solvable2TSP}. The results are summarised in tables and charts.

Summary section concludes the paper. Technical details and supportive material are placed in the appendices.

\section{Polynomially solvable case of the balanced  2TSP.}
\label{sec:solvable2TSP}

{\bf Balanced 2TSP with a Kalmanson distance matrix.}
\label{sec:recursions}
In this sub-section we consider the balanced 2TSP problem with an $n\times n$ Kalmanson distance matrix $(c_{ij})$, i.e. the matrix that satisfies  inequalities (\ref{kalm1.c})-(\ref{kalm2.c}).  An additional input for the problem is a subset $S$ of customers $S\subset\{1,2,\ldots,n\}$ to be visited in both periods. We will refer to these customers as \emph{fixed} nodes. For simplicity of notations we assume that $n+|S|$ is even and equals $2p$. The objective is to find two tours with the minimal total length such that the nodes from subset $S$ appear in each of the two tours, and the nodes from $\{1,\ldots,n\}\setminus S$ appear only in one of the tours. W.l.o.g.\ we assume that $1\in S$: we can assume it, since a cyclic renumbering of nodes does not change the property of a distance matrix being a Kalmanson matrix. The balancing condition demands that each of the two tours contains exactly $p$ nodes.

The tours $\seq{1,2,\ldots,n-1,n,1}$ and $\seq{1,n,n-1,\ldots,2,1}$ are the \emph{master} tours for the TSP with a Kalmanson matrix. It means that we can search for an optimal solutions of the balance 2TSP among the solutions where the nodes in each of the two tours are sorted either in increasing or decreasing order. What is to be decided is to how to allocate the nodes from $\{1,\ldots,n\}\setminus S$ into two tours.
We represent a feasible solution to the 2TSP as one sequence which is similar to a pyramidal tour. We start at node $1$, visit $p-1$ nodes from $\{ 2,\ldots,n\}$ including all nodes from $S$ in increasing order, arrive at node 1, then go from node 1 through nodes from $S$ and all nodes that have not yet been visited. In the second part of this sequence, i.e.\ after the second visit of node 1, we visit the nodes in decreasing order.

According to this description, the tours depicted in Fig.~\ref{fig:2pTSP} can be represented as the sequence
$\seq{1,2,3,5,8,9,10,{\bf 1},8,7,6,5,4,3,1}$.  It is important that  node ${\bf 1}$ is placed exactly in the middle of the sequence, with the nodes from $S$ placed before and then after ${\bf 1}$.

To construct any feasible sequence, we can start with a partially created sequence $\seq{1,\ldots,1}$ which contains only start and end node. There are following options for placing 2 in the sequence: (a) if node 2 is a fixed node, then we place one copy of $2$ after the start node $1$, and another copy of $2$ before the end node $1$: $\seq{1,2,\ldots,2,1}$; (b) if 2 is not a fixed node, we have two choices for the placement of 2 - either next to the start: $\seq{1,2,\ldots,1}$,  or next to the end: $\seq{1,\ldots,2,1}$. After placing 2, we decide on location of $3$, and so on. On each step of the construction, the feasible sequence contains two partially constructed tours for the 2TSP. First partially constructed tour starts from start node 1, the second partially constructed tour ends with the end node 1.

To find an optimal sequence among all feasible sequences, we use an approach similar to the one we used in finding an optimal pyramidal tour. Assume that we have already placed nodes $\{1,2,\ldots i,\ldots, j\}$, $i<j$. When we place $j+1$, 
 we have to consider two possibilities. If $j+1$ is not a fixed node, we place it in either first or second part of the sequence; the option with the shortest total length has to be chosen; if $j+1$ is a fixed node, then there is only one option to proceed - to place it in the both parts of the sequence. 
Assume that the second part of the partially constructed sequence (the second partially constructed tour)  contains $m, m\le p$, nodes. Let $V(i,j,m)$ be the minimal length of a feasible sequence that starts at node $i$, goes through the nodes $\{j+1,\ldots,n\}\cup\{1\}$ and stops at node $j$. Values $V(j,i,m)$ are defined similarly for the sub-tours that start at $j$ and end at $i$. For the case when $j$ is a fixed node, we also need to define $V(j,j,m)$ - the length of the shortest sequence from $j$ to $j$. 

It follows from the comments above that  values $V(i,j,m)$ and $V(j,i,m)$ satisfy the following dynamic programming equations:

\begin{equation}
\label{eq:vi1}
\begin{aligned}
V(i,j,m)=\min
\begin{cases}
  \min
  \begin{cases}
    c_{i,j+1}+V(j+1,j,m)\\
    c_{j+1,j} +V(i,j+1,m+1),\ \ \ \  \mbox{if}\ j+1\notin S;
     \end{cases}\\  
     c_{i,j+1}+c_{j+1,j}+V(j+1,j+1,m+1),\ \  \mbox{if}\ j+1 \in S.
\end{cases}\\
\qquad i=1,2,\ldots,j-1; j=i+1,\ldots,n-1.\\
V(j,i,m)=\min
\begin{cases}
  \min
  \begin{cases}
    c_{j,j+1}+V(j+1,i,m)\\
    c_{j+1,i} +V(j,j+1,m+1),\ \ \ \  \mbox{if}\  j+1\notin S;
     \end{cases}\\  
     c_{j,j+1}+c_{j+1,i}+V(j+1,j+1,m+1),\ \  \mbox{if}\ j+1 \in S.
\end{cases}\\
i=1,2,\ldots,j-1; j=i+1,\ldots,n-1.\\
\end{aligned}
\end{equation}
Boundary conditions (\ref{eq:vi2}) for the above recursions  define values that involve node $n$ and node~1. These conditions  ensure that the solution found is balanced:
\begin{equation}
\label{eq:vi2}
\begin{aligned}
V(i,n,m)=
  \begin{cases}
    c_{i1}+c_{1n} \ \ \mbox{if}\ m=p,\\
    \infty,\ \ \ \  \mbox{otherwise.}
 \end{cases}\\  
V(n,i,m)=
  \begin{cases}
    c_{n1}+c_{1i} \ \ \mbox{if}\ m=p,\\
    \infty,\ \ \ \  \mbox{otherwise.}
   \end{cases}\\ 
\qquad i=1,2,\ldots,n;\\
\end{aligned}
\end{equation}
Clearly, the length of an optimal sequence, i.e. the total length of an optimal pair of tours, can be calculated as 
\begin{equation}
\label{eq:vi3}
\begin{aligned}
L_{opt}=
  \begin{cases}
  min
  \begin{cases}
    c_{12}+V(2,1,1) \ \ \\
   c_{21}+V(1,2,2),\  \mbox{if}\ 2 \notin S;
 \end{cases}\\
 c_{12}+c_{21}+V(2,2,2),  \ \ \mbox{if}\ 2 \in S.
  \end{cases} 
  \end{aligned}
\end{equation}
This completes the proof of the main result of the section:
\begin{proposition}
An optimal balanced 2TSP problem with a Kalmanson distance matrix can be found in $O(n^3)$ time.
\end{proposition}

To solve the dynamic programming recursions (\ref{eq:vi1})-(\ref{eq:vi3}) one needs $O(n^3)$ space. The recursions can be re-written in a way that only $O(n^2)$ space is needed. To simplify reading of the paper, we have moved the corresponding recursions to Appendix A and leave the proof of the equivalence of the two systems as an exercise for an enthusiastic reader.

\bigskip
{\bf Generator of instances with known optimal solutions.} Easy solvable special cases of hard optimization problems can be used for generating test problems with known optimal solutions. Given a Kalmanson matrix, one can apply recursions (\ref{eq:vi1})-(\ref{eq:vi3}) and thus obtain an instance of the balanced 2TSP with the known optimal solution.  In this sub-section we describe an algorithm for generating Kalmanson matrices. To ensure the diversity of matrices generated, we  incorporate a random number generator in our procedure. 

As was shown in \cite{BSurv},  a symmetric $n\times n$ matrix $C$ is a Kalmanson matrix if and only if 
\begin{eqnarray}
c_{i,j+1}+c_{i+1,j} \le c_{ij}+c_{i+1,j+1}
&& \mbox{~for all~~}1\le i\le n-3,~ i+2\le j\le n-1 \label{kalmx.c1}  \\
c_{i,1}+c_{i+1,n} \le c_{in}+c_{i+1,1}
&& \mbox{~for all~~} 2\le i\le n-2. \label{kalmx.c2}
\end{eqnarray}
Let $\alpha_{ij}:=c_{ij}+c_{i+1,j+1}-c_{i,j+1}-c_{i+1,j}$ and
$\beta_i:=c_{in}+c_{i+1,1}-c_{i,1}-c_{i+1,n}$. Conditions (\ref{kalmx.c1})-(\ref{kalmx.c2}) are equivalent to the conditions 
$\alpha_{ij}>=0,\  \mbox{~for all~~}1\le i\le n-3,~ i+2\le j\le n-1$ and $\beta_i>=0, \mbox{~for all~~} 2\le i\le n-2$. Notice that if the values of any three out of four items involved in the definitions of $\alpha_{ij}$ and $\beta_i$ are known, then the value of the fourth item can always be defined to ensure the non-negativity of $\alpha_{ij}$ and $\beta_i$.


Given $(n-3)(n-2)/2$ non-negative numbers $\alpha_{ij}$ and $(n-3)$ non-negative numbers $\beta_i$, we show how to randomly generate Kalmanson matrices $C=(c_{ij})$ with 
$c_{ij}+c_{i+1,j+1}-c_{i,j+1}-c_{i+1,j}=\alpha_{ij}$ and $c_{in}+c_{i+1,1}-c_{i,1}-c_{i+1,n}=\beta_i$. The procedure will involve the following \emph{five} steps.

\begin{itemize}
\item{Step 1.} Generate randomly values $c_{1i}$, $i=2,\ldots,n$  and value $c_{2,n}$. Notice that the diagonal items are not involved in the definition of Kalmanson matrices, so we will always set them to $0$.
\end{itemize}

\begin{itemize}
\item{Step 2.} Define remaining items in the last column as
$c_{i,n}:=c_{i,1}+c_{i-1,n}-c_{i-1,1}-\beta_{i-1}$, for $i=3,\ldots,n-1$. Recall that $c_{1n}$ and $c_{2n}$ have already been defined on Step~1.
\end{itemize}

As an illustration to the procedure we consider an example with a $5 \times 5$ matrix. Assume that at Step 1 we have generated randomly $5$ integers: $3,1,3,0,2$, and obtained in such a way a partially constructed matrix $C^{(1)}$. Assume that $\beta_2=1$ and $\beta_3=1$, then Step 2 gives us matrix $C^{(2)}$, which has now first and last rows and first and last columns defined.

\[
C^{(1)}=
\begin{blockarray}{rrrrrr}
& 1 & 2 & 3 & 4 & 5 \\
\begin{block}{c(rrrrr)}
1\ \ \ & 0 & 3 & 1 & 3 & 0\\
2\ \ \ & 3 & 0 & \_ & \_ & 2 \\
3\ \ \ & 1 & \_ & 0 & \_ & \_ \\
4\ \ \ & 3 & \_ & \_ & 0 & \_ \\
5\ \ \ & 0 & 2 & \_ & \_ & 0 \\
\end{block}
\end{blockarray}\qquad\quad
C^{(2)}=
\begin{blockarray}{rrrrrr}
& 1 & 2 & 3 & 4 & 5 \\
\begin{block}{c(rrrrr)}
1\ \ \ & 0 & 3 & 1 & 3 & 0\\
2\ \ \ & 3 & 0 & \_ & \_ & 2 \\
3\ \ \ & 1 & \_ & 0 & \_ & -1 \\
4\ \ \ & 3 & \_ & \_ & 0 & 0 \\
5\ \ \ & 0 & 2 & -1 & 0 & 0 \\
\end{block}
\end{blockarray}\qquad\quad
\]

\begin{itemize}
\item{Step 3.} Define remaining items in the matrix, starting from row $i=2$; then define items in row $3$, and so on, until the single item in row $n-2$ is defined. In each row $i$, start with defining the item in column $j=n-1$, then define the item in column $n-2$, and so on, until the item in column $i+1$ is defined. The items are calculated as $c_{i,j}:=c_{i-1,j}+c_{i,j+1}-c_{i-1,j+1}-\alpha_{i-1,j}$.

\item{Step 4.}
If there are negative items in the matrix, find out the most negative item and subtract it from all items in the matrix (leaving diagonal items to remain $0$). This will result in a Kalmanson matrix with all  elements being non-negative.
\end{itemize}

For illustration, assume that
$\alpha_{13}=2$, $\alpha_{14}=3$, and $\alpha_{24}=1$. 
After defining all elements as prescribed at Steps 3, we obtain a symmetric Kalmanson matrix $C^{(3)}$ which is then transformed into a non-negative matrix $C^{(4)}$.
\[
C^{(3)}=
\begin{blockarray}{rrrrrr}
& 1 & 2 & 3 & 4 & 5 \\
\begin{block}{c(rrrrr)}
1\ \ \ & 0 & 3 & 1 & 3 & 0\\
2\ \ \ & 3 & 0 & -2 & 2 & 2 \\
3\ \ \ & 1 & -2 & 0 & -2 & -1 \\
4\ \ \ & 3 & 2 & -2 & 0 & 0 \\
5\ \ \ & 0 & 2 & -1 & 0 & 0 \\
\end{block}
\end{blockarray}\qquad\quad
C^{(4)}=
\begin{blockarray}{rrrrrr}
& 1 & 2 & 3 & 4 & 5 \\
\begin{block}{c(rrrrr)}
1\ \ \ & 0 & 5 & 3 & 5 & 2\\
2\ \ \ & 5 & 0 & 0 & 4 & 4 \\
3\ \ \ & 3 & 0 & 0 & 0 & 1 \\
4\ \ \ & 5 & 4 & 0 & 0 & 2 \\
5\ \ \ & 2 & 4 & 1 & 2 & 0 \\
\end{block}
\end{blockarray}\qquad\quad
\]

For test problems to be more challenging, it is possible to disguise the special structure of Kalmanson matrices by permuting rows and columns. In such a way one obtains the so-called \emph {permuted} Kalmanson matrix.
A matrix $A=(a_{ij})$ is called a permuted Kalmanson matrix if there exists a permutation $\sigma$ such that permuting rows and columns in $A$ yields the Kalmanson matrix $(a_{\sigma(i),\sigma(j)})$. To obtain a test instance with a random permuted Kalmanson matrix, we add the next step.

\begin{itemize}
\item{Step 5.} Generate a random permutation $\phi$ on the set $\{1,2,\ldots,n\}$ and permute rows and columns of the matrix to disguise the Kalmanson structure: the first row/column in the permuted matrix is now $\phi(1)$, the second row/column is $\phi(2)$, and so on.
\end{itemize}

For illustration, assume that a random generator produced a permutation $\phi=\{3,2,4,1,5\}$. After permuting rows and columns in Kalmanson matrix $C^{(4)}$ we obtained a \emph {permuted} Kalmanson matrix $C_{\phi}$, which does not satisfies conditions (\ref{kalmx.c1})-(\ref{kalmx.c2}) any more:

\[
C_{\phi}=
\begin{blockarray}{rrrrrr}
& 3 & 2 & 4 & 1 & 5 \\
\begin{block}{c(rrrrr)}
3\ \ \ & 0 & 0 & 0 & 3 & 1\\
2\ \ \ & 0 & 0 & 4 & 5 &4 \\
4\ \ \ & 0 & 4 & 0 & 5 & 2 \\
1\ \ \ & 3 & 5 & 5 & 0 & 2 \\
5\ \ \ & 1 & 4 & 2 & 2 & 0 \\
\end{block}
\end{blockarray}\qquad\quad
\]

Thus, to generate a random permuted Kalmanson matrix, one needs to generate non-negative numbers $\alpha_{ij}$ and $\beta_i$, and then to follow Steps 1-5 described above.

\bigskip
{\bf Heuristics motivated by the special solvable case.}
Given a symmetric $n\times n$ matrix $C$ and a subset $S$ of customers to be visited twice, we can find a \emph{feasible} solution to the 2TSP by using recursions (\ref{eq:vi1})-(\ref{eq:vi3}). This solution will be the best among exponentially many feasible solutions, however it not necessarily has to be the optimal solution. Even if matrix $C$ is a permuted Kalmanson matrix, one needs to permute $C$ into the Kalmanson matrix to guarantee that the solution found is the optimal solution. Clearly, permuting rows and columns in $C$ yields a new exponential neighbourhood searched in (\ref{eq:vi1})-(\ref{eq:vi3}), and hence a new feasible solution. Motivated by these simple observations, we suggest a heuristic approach to solving the 2TSP. The approach is based on the idea of enumerating special permutations with a hope that one of these permutations would permute $C$ into a matrix, which is ``close'' to a Kalmanson matrix.
   
  

It was noticed in \cite{ChFaTr}, that if a permuted Kalmanson matrix $C$ is transformed into matrix $C'$ as $c'_{i,j}=c_{i,j}-c_{i,1}-c_{1,j}$, i.e. into the matrix with zeros in the first row and column, then the $(n-1)\times(n-1)$ matrix, which is obtained by deleting the first row and first column in $C'$,  is a \emph{permuted} anti-Robinson matrix. 

A symmetric matrix $A=(a_{ij})$ is an  \emph{anti-Robinson matrix}, if for all $i<j<k$ it satisfies 
the conditions $a_{ik}\ge \max \{a_{ij},a_{jk}\}$, i.e.\ the entries in the matrix are placed in 
non-decreasing order in each row and column when moving away from the main diagonal (see \cite{PF} for more on anti-Robinson matrices and extensive references to the related publications).

The problem of recognising permuted Kalmanson matrices thus is reduced to recognising permuted anti-Robinson matrices (see also\cite{DRW} and  \cite{ChFaTr} for alternative recognition algorithms for Kalmanson matrices). Unfortunately all algorithms for recognising permuted anti-Robinson matrices or permuted Kalmanson matrices are quite elaborate (see e.g.\ \cite{PF}). We will use here a simple recognition algorithm that recognises only a sub-class of Kalmanson matrices, but has a reasonable intuitive motivation to be also used as a heuristic for the TSP.

We will call a Kalmanson matrix a \emph{strong} Kalmanson matrix, if all inequalities in (\ref{kalm1.c})-(\ref{kalm2.c}) are strict inequalities. Our algorithm will recognise permuted strong Kalmanson matrices; the algorithm looks very much like a nearest neighbour type heuristic for the TSP.

It follows from the definition of an anti-Robinson matrix, that the near-diagonal entries in the matrix correspond to the minimal Hamiltonian path in a complete graph with the weights of edges defined by the matrix. Moreover, this set of edges is also the minimal spanning tree for the graph (see \cite{ChFaTr}). So, a greedy algorithm for constructing a minimal Hamiltonian path as a spanning tree may lead to finding a permutation for permuting matrix into the anti-Robinson (Kalmanson) matrix. Motivated by this observation we suggest the following recognition-like algorithm.

The algorithm takes as an input an $n\times n$ matrix $C$, and returns a permutation $\tau$, which is a heuristic solution to the TSP with a distance matrix $C$. This permutation will also transform a permuted matrix into Kalmanson matrix, in cases when the inequalities (\ref{kalm1.c})-(\ref{kalm2.c}) are strict inequalities. So the algorithm can also be viewed as a recognition algorithm for permuted strong Kalmanson matrices. We will call our algorithm Kalmanson Nearest Neighbour (KNN).

\begin{tabbing}
1111\=1111\=1111\=111\=111\kill
\bf Algorithm KNN ($n$, $C$, $\tau$, start)\\
{\bf \{}
\>Transform $C$ into $C'$: $c'_{i,j}=c_{i,j}-c_{i,1}-c_{1,j}$, $i,j=2,\ldots,n$;\\
\>Initialise Hamiltonian path as $H=\{\mbox{\bf start}\}$;\\ 
\> Define current first and last nodes in $H$ as $\mbox{{\bf first}={\bf start}}$, $\mbox{{\bf last}={\bf start}}$;\\
\>\emph{ \bf repeat}  $n-2$ times\\
\>\emph{ \bf \{ }\\
\>\>Among nodes $\{2,3,\ldots,n\}$, which are not in $H$, find node $x$, which \\
\>\> is the nearest node to the set $\{\mbox{{\bf first}, {\bf last}}\}$; Use distance matrix $C'$;\\
\>\> \emph{{\bf if}} ($x$ is nearest to $\mbox{{\bf first}}$) \emph{{\bf then}}\\
\>\>\> $H=\seq{x,H}$, $\mbox{{\bf first}}=x$;\\
\>\>\emph{{\bf else}}\\
\>\>\>$H=\seq{H,x}$,  $\mbox{{\bf last}}=x$;\\
\>\emph{ \bf \} }\\
\>Define  $\tau=\seq{H}$; Insert $1$ into  $\tau$ at the position with the minimal \\ 
\>increase in the length of $\tau$ as a cycle. Use initial distance matrix $C$;\\
\> \emph{{\bf return }} tour $\tau$;\\
{\bf \}}
\end{tabbing}

\emph{Remark}. It can be proved that if matrix $C$ is a permuted \emph{strong} Kalmanson matrix, then there exists only one minimal Hamiltonian cycle, and hence $\tau$ is a unique permutation that permutes $C$ into the Kalmanson matrix $C_{\tau}$. If we were aiming only at a recognition of a Kalmanson matrix, the output of the algorithm could have been $\tau=\seq{1,H}$. So, the step of optimal insertion of $1$ could have been omitted. We know though that the permutation for permuting a matrix into a Kalmanson matrix is also an optimal tour for the TSP. So an insertion of $1$ into optimal position would not destroy the property of the permutation to recognise the matrix, but it adds more sense for the last step in the algorithm, if the initial matrix is not the Kalmanson one.   

\smallskip
We illustrate the algorithm on an instance with 5 nodes with the distance matrix $C_{\phi}$ used in the previous sub-section. The numbers  shown next to the matrix rows and columns are labels of the nodes.
We keep these labels as a reminder of the numbering of rows and columns in a Kalmanson matrix considered earlier.
To transform $C=C_{\phi}$ into $C'$ as prescribed in KNN algorithm, we subtract value $3$ from the forth row and column, subtract $1$ from fifth row and column, and obtain a matrix $C'$ with zeros in the first row and column. 

\[
C=C_{\phi}=
\begin{blockarray}{rrrrrr}
& 3 & 2 & 4 & 1 & 5 \\
\begin{block}{c(rrrrr)}
3\ \ \ & 0 & 0 & 0 & 3 & 1\\
2\ \ \ & 0 & 0 & 4 & 5 &4 \\
4\ \ \ & 0 & 4 & 0 & 5 & 2 \\
1\ \ \ & 3 & 5 & 5 & 0 & 2 \\
5\ \ \ & 1 & 4 & 2 & 2 & 0 \\
\end{block}
\end{blockarray}\qquad\quad
C'=
\begin{blockarray}{rrrrrr}
& 3 & 2 & 4 & 1 & 5 \\
\begin{block}{c(rrrrr)}
3\ \ \ & 0 & 0 & 0 & 0 & 0\\
2\ \ \ & 0 & 0 & 4 & 2 &3 \\
4\ \ \ & 0 & 4 & 0 & 2 & 1 \\
1\ \ \ &0 & 2 & 2 & 0 & -2 \\
5\ \ \ & 0 & 3 & 1 & -2 & 0 \\
\end{block}
\end{blockarray}\qquad\quad
\]

Assume that the start node chosen has label $2$, so $H=\seq{2}$.   The nearest node to $2$ has label $1$ (notice that label 3 is excluded from consideration). So, after adding next node to the Hamiltonian path we have $H=\seq{2,1}$ . The next node added to the path has label $5$: $H=\seq{2,1,5}$. The only left node is $4$, which is close to $5$, so $H=\seq{2,1,5,4}$. To decide on the best position for $3$ we have to find minimum among the following values (the indices are labels in $C$):
$c_{32}+c_{43}-c_{42}=0+0-4=-4$, $c_{23}+c_{31}-c_{21}=-2$, $c_{13}+c_{35}-c_{15}=2$, and $c_{53}+c_{34}-c_{54}=-1$. 
It is easy to see, there is only one position for $3$, and the output of the algorithm is $\tau=\seq{3,2,1,5,4}$, which is a permutation that permutes the matrix into the Kalmanson matrix. Recall that any of cyclic shifts of $\seq{1,2,3,4,5}$ or 
$\seq{5,4,3,2,1}$ permutes the matrix into the Kalmanson matrix.

\begin{figure}
\unitlength=1cm
\definecolor{qqqqff}{rgb}{0.,0.,1.}
{\begin{tikzpicture}[line cap=round,line join=round,>=triangle 45,x=1.0cm,y=1.0cm]
\clip(-4.,1.5) rectangle (2.34,4.68);
\draw (-3.,3.)-- (-1.,4.);
\draw (-1.,4.)-- (1.,3.);
\draw (1.,3.)-- (-1.,2.);
\draw (-3.,3.)-- (-1.,2.);
\begin{scriptsize}
\draw [fill=qqqqff] (-3.,3.) circle (2.5pt);
\draw[color=qqqqff] (-2.86,3.36) node {$1$};
\draw [fill=qqqqff] (-1.,4.) circle (2.5pt);
\draw[color=qqqqff] (-0.86,4.36) node {$2$};
\draw [fill=qqqqff] (-1.,2.) circle (2.5pt);
\draw[color=qqqqff] (-0.86,2.36) node {$4$};
\draw [fill=qqqqff] (1.,3.) circle (2.5pt);
\draw[color=qqqqff] (1.14,3.36) node {$3$};
\end{scriptsize}
\end{tikzpicture}
}
\put(6.5,1.5){\makebox(3,0)[cc]
{
\[
D=(d_{ij})=
\begin{blockarray}{rrrr}
\begin{block}{(rrrr)}
 0 & 22 & 40 & 22 \\
22 & 0 & 22 & 20  \\
 40 & 22 & 0 & 22 \\
 22 & 20 & 22 & 0  \\
\end{block}
\end{blockarray}\qquad\quad
\]
}}
\caption{Illustration for the KNN-heuristic: 4 points on the convex hull and the distance matrix.}
\label{fig:diamand}
\end{figure}
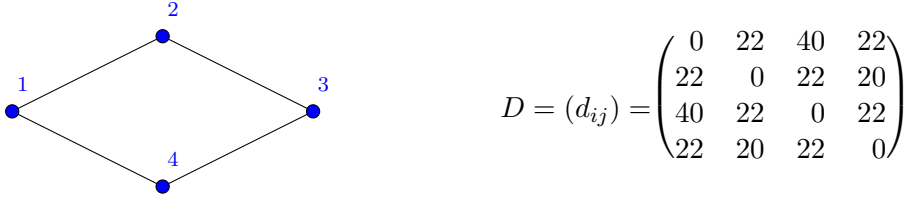

\smallskip
To illustrate an interesting feature of the KNN-algorithm as the TSP-heuristic, we refer to the instance shown in Fig.~\ref{fig:diamand}. If one would apply the classic nearest neighbour TSP heuristic to the instance, this heuristic will always include edge $(2,4)$ in the TSP tour, independently on the starting node. Contrary to this, the KNN-heuristic will always find the tour $\seq{1,2,3,4,1}$ (or its inverse and/or cyclic shift equivalents).
Indeed, after transforming matrix $D$ into matrix $D'$ with zeros in the first row and column, we have $d'_{23}=d'_{32}=-40$, $d'_{24}=d'_{42}=-24$, and
$d'_{34}=d'_{43}=-40$.  If we choose, for instance, $3$ as the start node, then the Hamiltonian cycle will grow as $H=\{3\}$, $H\{3,2\}$, and $H=\{4,3,2\}$. 
To insert $1$ into the best position we have to compare $d_{21}+d_{14}-d_{24}=22+22-20=24$, $d_{41}+d_{13}-d_{43}=22+40-22=40$, and $d_{31}+d_{12}-d_{32}=40+22-22=40$. So the best insertion would lead to the tour $\seq{1,4,3,2,1}$, which also defines the locations of the points on their convex hull.

\smallskip
We are ready to describe now a simple heuristic for solving the balanced 2TSP with an arbitrary distance matrix.
We will call this heuristic Kalmanson Sequence (KS) Heuristic.

\begin{tabbing}
1111\=1111\=1111\=111\=111\kill
\bf KS-Heuristic for the balanced 2TSP\\
{\bf \{}
\>\emph{ \bf for}  \emph{start}$=2,\ldots,n$ { \bf do} \\
\>\emph{ \bf \{ }\\
\>\> Apply KNN($n$, $C$, $\tau$, \emph{start}) to find permutation $\tau$;\\
\>\>Construct two sub-tours by applying recursions (\ref{eq:vi1})-(\ref{eq:vi3}) to the permuted matrix $C_{\tau}$;\\
\>\> Apply Belperm TSP heuristic for each of the sub-tours; Save the record;\\
\>\emph{ \bf \} }\\
\>{\bf return} record - the best solution found in $n-1$ iterations;\\
{\bf \}}
\end{tabbing}



\bigskip
{\bf Computational experiments with KS-heuristic.}\label{sec:experiments1}
In this sub-section we report results of our empirical testing of the KS-heuristic. The only published benchmark instances for the balanced 2TSP are instances from Bassetto \& Mason \cite{BassettoM11} paper. The instances with Kalmanson matrices (i.e.\ with known optimal solutions) are not useful here, since the KS-Heuristic was designed to find optimal solutions for this type of instances.

\smallskip

 Bassetto \& Mason \cite{BassettoM11}  considered a Euclidean version of the balanced 2TSP.
A short summary of the approaches used in \cite{BassettoM11} is as follows. 

First a TSP tour on the set of all customers, called a general tour (GT), is constructed. The GT is used to obtain a partition of customers into two subsets visited in two periods. The initial partition of GT into two sub-tours is improved by applying decision rules motivated by geometry (e.g. solutions have been improved by removing crossing edges).
For each of the sub-tours, an optimal TSP tour is constructed by applying an exact TSP algorithm (see \cite{ABCC}, Chapter 16,  and \cite{Concorde}). 
The authors also used visualisation and human intervention for improvements of the solutions found by a computer. 

The set of benchmark instances from \cite{BassettoM11} contains 60 randomly generated instances with 48 customers. The set of these instances is divided into three subsets with a different number of customers to be visited in the two periods: $8$, $12$, and $24$ customers.   
We used these instances to test the KS-heuristic.
Summary of our computational experiments is presented in Table \ref{table:summaryHeuristics}.
We compare our results with the both types of solutions presented in \cite{BassettoM11}: best solutions found by a computer, which we label as ``PC" solutions, and  the solutions obtained after visualisation and manual improvements,  labelled as ``PC+manual" (or ``PC+m" in some tables). 
  For each set of 20 instances, Table~\ref{table:summaryHeuristics} shows the mean of the percentage above the best known solution, the percentages for the ``best" and ``worst" instance, and the number of solutions where we improved the best known result or found the same solution (``improved \#"). Negative percentage means that the best known solution was improved. Results for individual test instances can be found in  appendix B.

For our experiments we used a desktop computer with Intel i7-3770 3.40 GHz CPU, 16 GB of RAM, and GNU C++ compiler. Computational time for the instances was just a fraction of a second (``few seconds" computational time was mentioned in  \cite{BassettoM11}).

The results show that  the KS-heuristic, which is based on polynomially solvable case, is competitive with the previously published heuristics: for 44 out of 60 instances we found the same or better solution, for all three sets of instances the mean percentage  shows the better overall performance if compared with solutions found by computer. 

\begin{table}
\begin{center}
\begin{tabular}{||c||l||c|c||}
\hline 
\cline{2-4} 
& 
& PC & PC+manual  
\\ 
\hline 
&Mean \%   & -0.63  & 1.10  \\ 
8 nodes& Best \% & -5.45  & -0.91   \\ 
visited twice&Worst \%  &2.11  & 3.83   \\ 
&Improved \# & 14/20 & 2/20 \\ 
\hline 
&Mean \%  & -0.42  & 0.62  \\ 
16 nodes& Best \%  & -2.71  & -1.03   \\ 
visited twice&Worst \%  & 3.14  & 3.14   \\ 
&Improved \# & 15/20 & 7/20 \\ 
\hline 
&Mean \%  & -0.32  & 0.30  \\ 
24 nodes& Best \%   & -1.8  & -1.42   \\ 
visited twice&Worst \%   & 1.13  & 1.91   \\ 
&Improved \# & 15/20 & 7/20 \\ 
\hline 
Total&Improved \#  & 44/60 & 16/60 \\ 
\hline 
\end{tabular} 
\end{center}
\caption{Summary of computational experiments for the KS-heuristic.}
\label{table:summaryHeuristics}
\end{table}

\bigskip
In this section we considered a version of the 2VRP known as the balanced 2TSP. We described a new  polynomially solvable case of the 2TSP and a new heuristic.
The key feature of our heuristic is a dynamic programming algorithm for solving recursions (\ref{eq:vi1})-(\ref{eq:vi3}). In the next section we describe an approach which is also based on the dynamic programming paradigm but can be used for solving 2VRP with much more complicated set of additional constraints.

\section{A framework for the 2VRP}
\label{sec:framework}
{ \bf  Notations and definitions.}
In this section we model the 2VRP in a way which allows us to easily incorporate a variety of additional constraints that can be found in real life applications.

We assume that there are two {\em heterogeneous} vehicles, vehicle 1 and vehicle 2. They have different capacities $W_1$ and $W_2$ and different mileage costs. The vehicles are used for delivering  goods to customers. The vehicles may travel from different depots and return to different depots. Assume that vehicle 1 travels from depot $d^1_R$ and returns to depot $d^1_L$, and vehicle 2 travels from depot $d^2_R$ and returns to depot $d^2_L$.

A set of $n$ customers $\{1,2,\ldots,n\}$, with the demands for delivery $w(i)$, $i=1,\ldots,n$, is given. The customers have to be visited by one of the two vehicles that delivers the demanded goods from a depot to the customers. The total demand of all customers is less than the total capacity of the two vehicles and can be partitioned in a way that only one route per vehicle is needed. Both depots have unlimited inventory.

We assume that each customer is located in an estate with only two entry points. A network of one-way roads within the estate connects these two points. So the travel costs within the estate are asymmetric costs.
To distinguish between two entry points to the estate, we refer to  one of the entry points as the {\em left node} and denote it by $L(i)$. The other of the two points is  referred as the {\em right node}, and is denoted by $R(i)$. 

Representation of locations in the described way makes the task of data aggregating much easier. Imagine that you have a collection of customers located on a street. If you envisage that these customers can be serviced by one vehicle, you substitute all these customers with one customer which has the demand equals to the total demand of the customers on the street. The location of the new customer $i$ is defined by two nodes, $L(i)$ and $R(i)$, which are the beginning and the end of the street.  The costs of travelling from $L(i)$ to $R(i)$, and from $R(i)$ to $L(i)$,  can be different. If the street is a one-way street, the one of the costs is to be set to the infinity.

With each customer $i$ we associate a set of seven attributes
$\{L(i)$, $R(i)$, $l^1_L(i)$, $ l^1_R(i)$, $ l^2_L(i)$, $l^2_R(i)$, $ w(i)\}$ with the following meaning:
\begin{itemize}
\item the left node $L(i)$;
\item the right node $R(i)$;
\item cost of travelling from the left node to the right node $l^m_L(i)$, if travelled by vehicle $m$, $m=1,2$;
\item cost of travelling from the right node to the left node $l^m_R(i)$, if travelled by vehicle $m$, $m=1,2$;
\item the demand $w(i)$.  
\end{itemize}

The costs of travelling between $n$ customers and 4 depots are given by the two $2(n+2)\times 2(n+2)$ cost matrices $C^1=(c^1_{ij})$ and $C^2=(c^2_{ij})$,  where $c^m(i,j)$ is the cost of travelling by vehicle $m$ from node $i$ to node $j$. 

A visit of customer $i$ adds to the travel costs either cost $l^m_L(i)$, or cost $l^m_R(i)$ depending on the way the customer is visited - from the left entry node, or from the right entry node.

The main idea of our approach is to view two routes for the two vehicles as \emph{one two-vehicle} route, and  utilise the well known dynamic programming algorithm for \emph{the TSP} for finding the optimal route.
To implement this idea, we introduce an auxiliary customer $0$ with the set of attributes $\{d^1_L,d^2_R,0,\infty,0,\infty,0\}$. The reason for introducing the auxiliary customer $0$ and placing it into the two-vehicle route is to separate points visited by the two different vehicles. It costs nothing to travel from $d^1_L$ to $d^2_R$ (therefore $l^1_L(0)=l^2_L(0)=0$), while travelling in the opposite direction is forbidden (by the infinitely large cost). The demand of this customer is obviously zero.

The two-vehicle route starts from node $d^1_R$, visits all customers  from the set $U=\{0\}\cup \{1,2,\ldots,n\}$, and ends in the node $d^2_L$. Vehicle 1 is the first vehicle visiting customers in the route. After visiting customer $0$ the mode of travelling is changed from travelling by vehicle 1 to travelling by vehicle 2. Vehicle 2 will travel from customer $0$, i.e.\ from the depot $d^2_R$ to the depot $d^2_L$ (recall that $R(0)=d^2_R$).

The objective is to find a minimum cost route  of delivering the requested demand to all customers. The total demand of customers visited by each vehicle should not exceed the corresponding vehicle capacity.

\smallskip
{\bf Dynamic programming recursions.}
In this section, the well known Held \& Karp \cite{HeldKarp} dynamic programming algorithm for the TSP is adapted for the case of the 2VRP model formulated above.

\begin{figure}
\unitlength=0.9cm
\definecolor{qqqqff}{rgb}{0.,0.,1.}
\begin{tikzpicture}[line cap=round,line join=round,>=triangle 45,x=1.0cm,y=1.0cm]
\clip(1.55,1.11) rectangle (17.,10.3);
\draw (2.,6.)-- (4.,6.);
\draw (6.007042253521121,6.535211267605635)-- (8.00704225352112,6.535211267605635);
\draw (6.,5.492957746478873)-- (8.,5.492957746478873);
\draw [rotate around={-3.576334374995288:(7.,6.001760563380284)},dash pattern=on 6pt off 6pt] (7.,6.001760563380284) ellipse (1.5307841759435992cm and 1.5179812710839782cm);
\draw [rotate around={-1.2364076028561437:(12.163732394366187,6.32218309859155)},dash pattern=on 6pt off 6pt] (12.163732394366187,6.32218309859155) ellipse (2.971769702905846cm and 2.428082112076715cm);
\draw [rotate around={-0.1469119331988365:(10.943661971830982,6.026408450704226)}] (10.943661971830982,6.026408450704226) ellipse (5.81427438136858cm and 3.271811342548349cm);
\draw (9.4,2.7) node[anchor=north west] {Set of customers $J$};
\draw (10.278169014084503,5.545774647887325) node[anchor=north west] {Set of customers $J \setminus \{j\}$};
\draw (6.,5.4) node[anchor=north west] {customer $j$};
\draw (11.95422535211267,7.616197183098592)-- (13.137323943661965,7.616197183098592);
\draw (11.239436619718305,8.330985915492958) node[anchor=north west] {$d^1_L$};
\draw (13.383802816901403,8.28169014084507) node[anchor=north west] {$d^2_R$};
\draw (11.63,7.4) node[anchor=north west] {customer $\mathbf{0}$};
\draw (3.9,7.3) node[anchor=north west] {$c^1_{R(i)L(j)}$};
\draw (3.9,5.6) node[anchor=north west] {$c^1_{R(i)R(j)}$};
\draw [dash pattern=on 1pt off 3pt on 6pt off 4pt] (4.,6.)-- (6.007042253521121,6.535211267605635);
\draw [dash pattern=on 1pt off 3pt on 6pt off 4pt] (5.128553288839242,6.300947543690467) -- (5.046389296616109,6.106849996844511);
\draw [dash pattern=on 1pt off 3pt on 6pt off 4pt] (5.128553288839242,6.300947543690467) -- (4.960652956905013,6.428361270761123);
\draw [dash pattern=on 1pt off 3pt on 6pt off 4pt] (4.,6.)-- (6.,5.492957746478873);
\draw [dash pattern=on 1pt off 3pt on 6pt off 4pt] (5.125433196241976,5.714678907994991) -- (4.959114330399999,5.585207620928324);
\draw [dash pattern=on 1pt off 3pt on 6pt off 4pt] (5.125433196241976,5.714678907994991) -- (5.0408856696,5.907750125550549);
\draw [dash pattern=on 1pt off 3pt on 6pt off 4pt] (8.00704225352112,6.535211267605635)-- (10.031690140845065,6.531690140845072);
\draw [dash pattern=on 1pt off 3pt on 6pt off 4pt] (9.14876740994213,6.533225658637946) -- (9.019076852856427,6.367077716392306);
\draw [dash pattern=on 1pt off 3pt on 6pt off 4pt] (9.14876740994213,6.533225658637946) -- (9.01965554150976,6.699823692058401);
\draw [dash pattern=on 1pt off 3pt on 6pt off 4pt] (8.,5.492957746478873)-- (9.982394366197179,5.471830985915494);
\draw [dash pattern=on 1pt off 3pt on 6pt off 4pt] (9.120591243742526,5.481015388641475) -- (8.989424211955537,5.316030573940693);
\draw [dash pattern=on 1pt off 3pt on 6pt off 4pt] (9.120591243742526,5.481015388641475) -- (8.992970154241641,5.648758158453672);
\draw (7.1,8.3) node[anchor=north west] {$+VL[j,J\setminus\{j\}]$};
\draw (7.3,4.7) node[anchor=north west] {$+VR[j,J\setminus\{j\}]$};
\draw (6.063380281690138,7.369718309859155) node[anchor=north west] {$L(j)$};
\draw (7.542253521126757,7.345070422535212) node[anchor=north west] {$R(j)$};
\draw (5.94014084507042,6.334507042253522) node[anchor=north west] {$R(j)$};
\draw (7.616197183098588,6.285211267605635) node[anchor=north west] {$L(j)$};
\draw (2.,6.) node[anchor=north west] {customer $i$};
\draw (1.57,6.7) node[anchor=north west] {$L(i)$};
\draw (2.637323943661972,7.073943661971831) node[anchor=north west] {$l^1_L(i)$};
\begin{scriptsize}
\draw [fill=qqqqff] (2.,6.) circle (2.5pt);
\draw [fill=qqqqff] (4.,6.) ++(-2.5pt,0 pt) -- ++(2.5pt,2.5pt)--++(2.5pt,-2.5pt)--++(-2.5pt,-2.5pt)--++(-2.5pt,2.5pt);
\draw [fill=qqqqff] (6.007042253521121,6.535211267605635) circle (2.5pt);
\draw [fill=qqqqff] (8.00704225352112,6.535211267605635) ++(-2.5pt,0 pt) -- ++(2.5pt,2.5pt)--++(2.5pt,-2.5pt)--++(-2.5pt,-2.5pt)--++(-2.5pt,2.5pt);
\draw [fill=qqqqff] (6.,5.492957746478873) ++(-2.5pt,0 pt) -- ++(2.5pt,2.5pt)--++(2.5pt,-2.5pt)--++(-2.5pt,-2.5pt)--++(-2.5pt,2.5pt);
\draw [fill=qqqqff] (8.,5.492957746478873) circle (2.5pt);
\draw [fill=qqqqff] (11.95422535211267,7.616197183098592) circle (2.5pt);
\draw [fill=qqqqff] (13.137323943661965,7.616197183098592) ++(-2.5pt,0 pt) -- ++(2.5pt,2.5pt)--++(2.5pt,-2.5pt)--++(-2.5pt,-2.5pt)--++(-2.5pt,2.5pt);
\end{scriptsize}
\end{tikzpicture}
\vspace{-0.4cm}
\caption{Illustration to calculation of $\VL[i,J]$, which is the length of the shortest route from customer $i$ through all customers in set $J$ (the case of $\mathbf{0}\in J$).}
\end{figure}

Let $J$ be a subset of customers not containing $i$, so $J\subset U$, $i \notin J$.
Let $\VL[i,J]$ be the minimum cost of an optimal 2-vehicle route among all routes that start visiting customer $i$ from the left node, then visiting all the customers in set $J$, and stopping in depot $d^2_L$. Similarly, define $\VR[i,J]$ to be the cost of the optimal route that starts visiting customer $i$ from the right node. The optimal cost of the 2-vehicle tour can be calculated as
 
\begin{equation}
\label{eq:dp0}
V=min_{i\in U\setminus \{0\}}\{c^1_{d^1_R,L(i)}+\VL[i,U\setminus \{i\}],c^1_{d^1_R,R(i)}+\VR[i,U\setminus \{i\}]\}.
\end{equation}

We assume here that the total demand of customers is bigger than the capacity of each of the vehicles, so in the formula above $0$ cannot be the first customer. 
Values $\VL[i,J]$ and $\VR[i,J]$ for all customers $i$ and subsets $J$, $J\subset \{0\}\cup\{1,2,\ldots,n\}$, are calculated as shown in the recursions below:

\begin{equation}
\label{eq:dp1}
\begin{aligned}
\VL[i,J]& =
\begin{cases}
\left.
\min_{j\in J}
\begin{cases}
				l^1_L(i)+c^1_{R(i),L(j)}+\VL[j,J\setminus\{j\}] \\
				l^1_L(i)+c^1_{R(i),R(j)}+\VR[j,J\setminus\{j\}]
\end{cases}
\hspace{-2.0ex}\right\}
&\mbox{if }0\in J, \\
\\
\left.
\min_{j\in J}
\begin{cases}
				l^2_L(i)+c^2_{R(i),L(j)}+\VL[j,J\setminus\{j\}]\\
				l^2_L(i)+c^2_{R(i),R(j)}+\VR[j,J\setminus\{j\}]
\end{cases}
\hspace{-2.0ex}\right\}
&\mbox{if }0\notin J, w(\{i\}\cup J)\le W_2\\
\infty & \mbox{otherwise}
\end{cases}
\end{aligned}
\end{equation}

\begin{equation}
\begin{aligned}
\VR[i,J]& =
\begin{cases}
\left.
\min_{j\in J}
\begin{cases}
				l^1_R(i)+c^1_{L(i),L(j)}+\VL[j,J\setminus\{j\}] \\
				l^1_R(i)+c^1_{L(i),R(j)}+\VR[j,J\setminus\{j\}]
\end{cases}
\hspace{-2.0ex}\right\}
&\mbox{if }0\in J, \\\\
\left.
\min_{j\in J}
\begin{cases}
				l^2_R(i)+c^2_{L(i),L(j)}+\VL[j,J\setminus\{j\}]\\
				l^2_R(i)+c^2_{L(i),R(j)}+\VR[j,J\setminus\{j\}]
\end{cases}
\hspace{-2.0ex}\right\}
&\mbox{if }0\notin J, w(\{i\}\cup J)\le W_2\\
\infty &\mbox{otherwise}
\end{cases}
\end{aligned}
\end{equation}

\begin{equation}
\begin{aligned}
\VL[0,J]& =
\begin{cases}
\left.
\min_{j\in J}
\begin{cases}
				c^2_{R(0),L(j)}+\VL[j,J\setminus\{j\}\\
				c^2_{R(0),R(j)}+\VR[j,J\setminus\{j\}
\end{cases}
\hspace{-2.0ex}\right\}
&\mbox{if }
\begin{cases}w(U\setminus (\{i\}\cup J))\le W_1,\\ 
\mbox{and }w(\{i\}\cup J)\le W_2
\end{cases}\\
\infty & \mbox{otherwise}
\end{cases}
\end{aligned}
\end{equation}

The boundary conditions are:
\begin{equation}
\label{eq:dpB}
\begin{aligned}
\VL[i,\emptyset]& =
				l^2_R(i)+c^2_{L(i),d^2_L},
\\
\VR[i,\emptyset]& =
				l^2_L(i)+c^2_{R(i),d^2_L}.
\end{aligned}
\end{equation}

Recursions (\ref{eq:dp0})-(\ref{eq:dpB}) extend Held \& Karp recursions to the case of the 2VRP. Since there are only two vehicles, the capacity constraints are easily verified without any extra dimension or complicated calculations. Notice that we use notation $w(J)$ for the sum of demands of all items in set $J$.

\bigskip
{\bf Aggregation strategy and local search.}
\begin{figure}
\unitlength=1cm
\begin{center}
\begin{picture}(11.5,7.5)
{
\begin{picture}(10.5,6.5)
\includegraphics[scale=0.8]{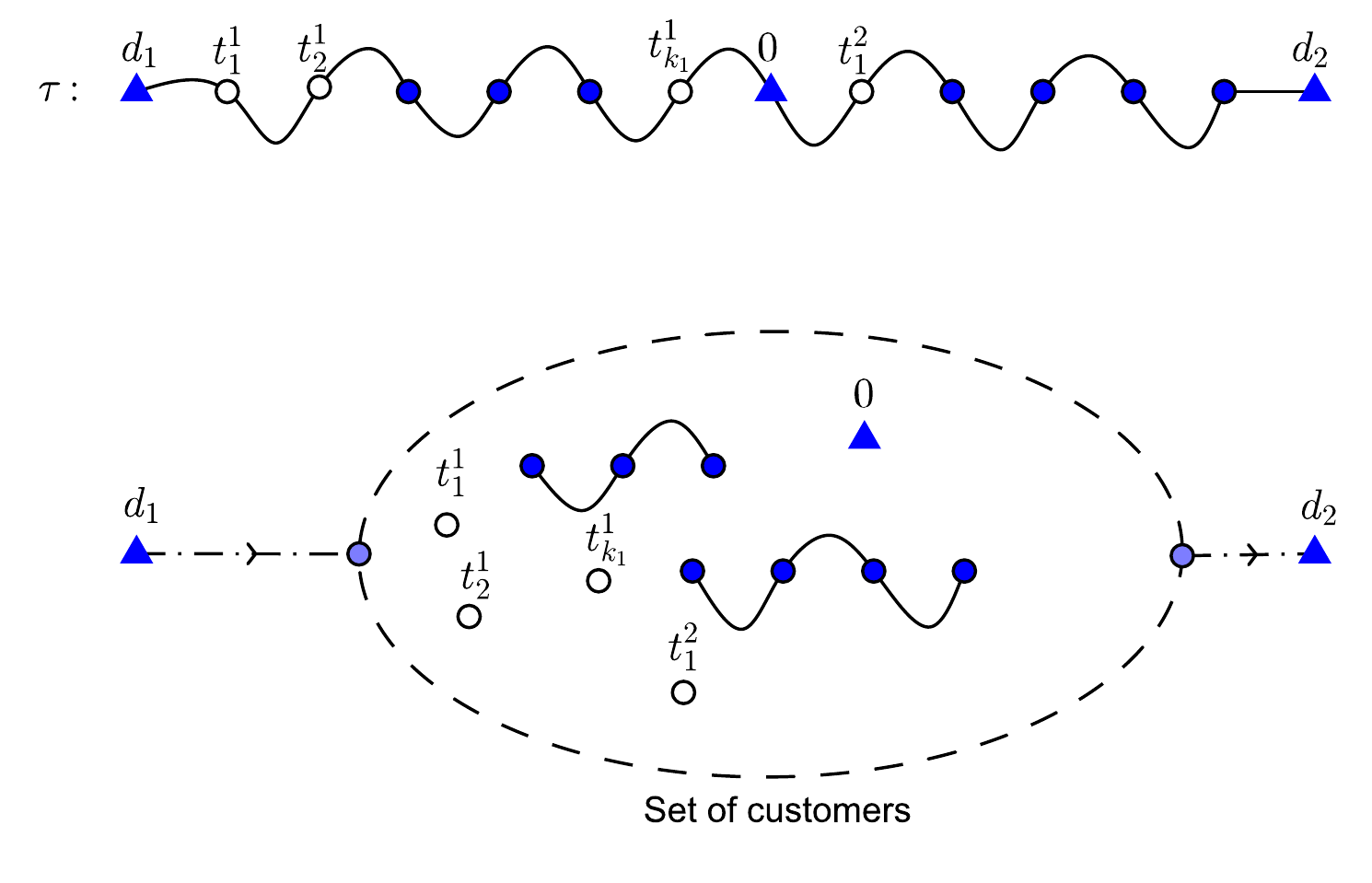}
\end{picture} 
}
\end{picture}
\end{center}
\caption{ Illustration of the ``sliding subsets" heuristic: first step of disassembling; each subset contains $2$ customers, $S_1=\{t_1^1,t_2^1\}$, $S_2=\{t_{k_1}^1,t_1^2\}$.}
\label{fig:disassembling1}
\end{figure}
The dynamic programming approach above has an exponential time complexity and therefore can be used only for small size problems.
In this section we describe an approach for data aggregation that can be used for reducing the size of the problems. 

Consider the 2VRP with $n$ customers, where $n$ is big enough to make the dynamic programming recursions computationally intractable. In this case we use the following approach.
First, we use a simple heuristic to find a feasible 2VRP solution.  We ``disassemble" then this solution into a small number of subpaths, and  
represent each subpath as a customer. 
Recalculation of the attributes associated with the new customer $i$ is straightforward. The left node $L(i)$ for the new customer is the first node  in the subpath (i.e. either left or right node of the first customer in the subpath); the right node $R(i)$ of the new customer is the last node in the subpath; the demand $w(i)$  is the sum of the demands of customers in the subpath; values $l_L^1(i)$, $l_L^2(i)$, $l_R^1(i)$, $l_R^2(i)$ are calculated as the travelling costs through the subpath in different directions and by different vehicles.

We apply then recursions (\ref{eq:dp0})-(\ref{eq:dpB}) and find an optimal solution to the new small-size 2VRP.   Obviously, an exact solution obtained can be viewed only as an approximate solution for the initial 2VRP. This solution can be  ``disassembled" again to get the new small-size problem. The process of solving small size problems is repeated until  all possible ways of disassembling and aggregating have been enumerated and no further improvement was achieved. 

\smallskip
We suggest the following straightforward approach for disassembling, which we call the sliding sub-set method. Assume that we have an initial solution to the 2VRP: $\tau=\seq{d_1,t^1_1,t^1_2,\ldots,t^1_{k_1},0,t^2_1,t^2_2,\ldots,t^2_{k_2},d_2}$. Here the route of vehicle 1 is 
$\seq{d_1,t^1_1,t^1_2,\ldots,t^1_{k_1},0}$ and the route of vehicle 2 is $\seq{0,t^2_1,t^2_2,\ldots,t^2_{k_2},d_2}$.
We disassemble this solution into a new set of customers as follows.

First, we leave customer $0$ to be the customer in the new 2VRP, and delete it from $\tau$.
Let $s$ now be a small constant, a parameter of the algorithm. Choose two subsets of customers $S_1$ and $S_2$  containing $s$ items each. The customers in each subset are on consecutive positions in $\tau$. Subset $S_1$ will always contain at least one customer from the route of vehicle 1, and $S_2$ contains at least one customer from the route of vehicle 2. 

On the first disassembling step  define  $S_1=\{t^1_1,\ldots,t^1_s\}$, and  $S_2=\{t^1_{k_1-s+2},\ldots,t^1_{k_1},t^2_1\}$. Notice that subset $S_2$ is chosen to ensure that at least one customer from vehicle 2 is included in the subset. Delete $S_1$ and $S_2$ from $\tau$ and add them to the set of customers in the new 2VRP. Leave $d_2$ to be the depot in the new problem, delete it from $\tau$. 

Sub-paths which are left in $\tau$ are to be replaced by  aggregated customers and to be added to the new 2VRP. On the first step these customers are: depot $d_1$, the subpath $\seq{t^1_{s+1},\ldots,t^1_{k_1-s+1}}$, and the subpath
$\seq{t^2_2,\ldots,t^2_{k_2}}$. The new 2VRP contains $2s +5$ customers (depots are counted as customers). 

The small size 2VRP obtained from the disassembled solution is solved to optimality. If no better solution is found, we \emph{redefine} subset $S_2$ by deleting, say, $l$ first elements from $S_2$ and adding $l$ new elements; $l$ is a parameter, to which we refer as {\em step}.  We repeat the process until we reach the end of the tour. 
When all feasible subsets $S_2$ are enumerated, we change set $S_1$ (with the step $l$) and redefine set $S_2$ to follow set $S_1$ similar to the first step described above.

If the solution is improved, the process of disassembling is applied to the new solution. The procedure stops when all feasible subsets $S_1$ and $S_2$ are enumerated and no improvements found.

Fig.~\ref{fig:disassembling1} illustrates the concept of sliding subsets. 
To simplify drawings, the customers in the initial 2VRP are depicted as points, not as intervals.

\begin{figure}
\unitlength=0.9cm
\begin{center}
{
\begin{picture}(12,8.0)
\includegraphics[scale=0.8]{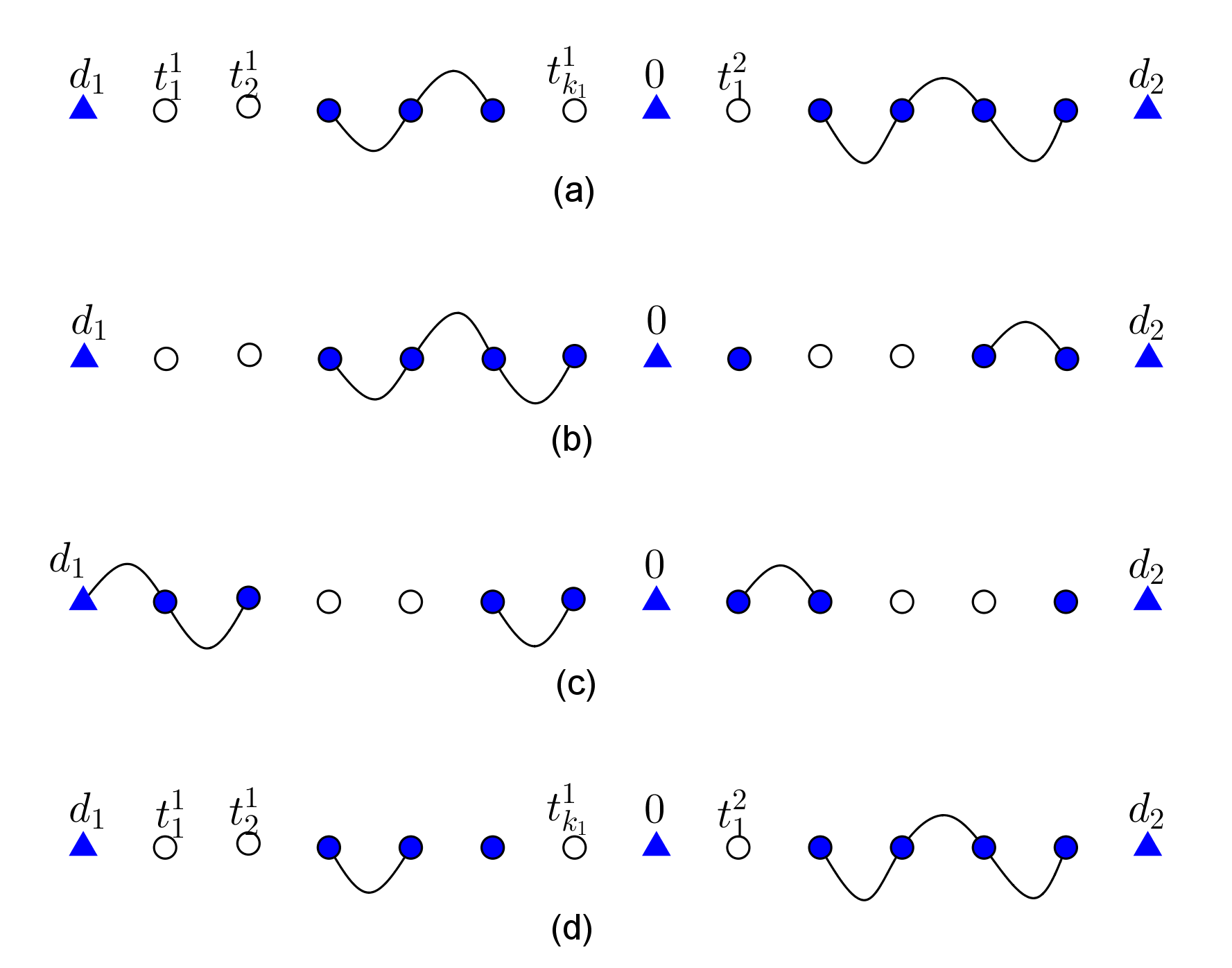}
\end{picture} 
}
\end{center}
\caption{ Illustrations of disassembling: (a) first step: $S_1$ and $S_2$ are separated by one sub-path only; (b) an example when $S_1$ and $S_2$ are separated by two sub-paths and a depot; (c) depot and the first subpath are considered as one customer; (d) modified first step: a subpath between $S_1$ and $S_2$ (case (a)) is partitioned into a sub-path and a single customer to keep the size of the new 2VRP fixed.
}
\label{fig:disassembling2}
\end{figure}

\bigskip
With the fixed parameter $s$, the disassembling procedure as described above would yield the small 2VRPs of different size.
To keep the size of the 2VRPs fixed and equal to $2s+6$, we slightly modify the procedure.
We refer to Fig.~\ref{fig:disassembling2} in our explanations of various steps (and possibilities) of the process, with the parameters $s=2$ and $l=2$.

 Fig.~\ref{fig:disassembling2}(a) illustrates the first step of the disassembling process as was described above. Sets $S_1$ and $S_2$ are separated by one sub-path in this case. The size of the small 2VRP  is $2s+5$.

Fig.~\ref{fig:disassembling2}(b) illustrates the outcome of disassembling the tour on the next iteration. Notice that depot $0$ is always treated as a separate customer, therefore setting $l=2$ yields the position of $S_2$ as shown in the figure.  There are two sub-paths between sets $S_1$ and $S_2$, and the size of the small 2VRP  is $2s+6$. 

Consider the step when $S_1$ and $S_2$ are chosen as shown in Fig.~\ref{fig:disassembling2}(c). If the first sub-path did not contain the depot, the size of the problem would have been $2s+7$.

So, with fixed parameter $s$ we may have 2VRPs with $2s+5$, $2s+6$, and $2s+7$ customers. On the implementation step, it was convenient to keep the size of the problem fixed at $2s+6$. Therefore it was decided (1) to ``glue" the first sub-path with the depot and define it as the depot in the new problem, as shown in Fig.~\ref{fig:disassembling2}(c); 
(2) in case when subsets  $S_1$ and $S_2$ are separated by a single path (for example, on the first step of disassembling), it was decided to consider the last node in the sub-path as a sub-path with one node: in this case the problem with $2s+5$ customers becomes the problem with $2s+6$ customers (compare Fig.~\ref{fig:disassembling2}(a) and Fig.~\ref{fig:disassembling2}(d)). 



\smallskip
{\bf Variety of the VRPs covered by the framework.} Below we illustrate the  advantages of the suggested dynamic programming approach by listing some types of the VRPs that can be tackled by the approach proposed.
\begin{itemize}
\item \emph{Arc routing.} Our model made no difference between the classic capacitated VRP, where the customers are represented by one node in a road network, and the arc routing VRP, where the customers are the streets/arcs in a network
(see \cite{Wohlk,CP}).

\item \emph{Heterogeneous fleet.} The dynamic programming recursions take into account individual characteristics of the vehicles, so both homogeneous and heterogeneous fleets (see \cite{BaldachiVigo} and recent survey \cite{KBJL})  can be managed.

\item \emph{Multi-depot and open VRPs.} Incorporation of the multi-depot feature  into the model is straightforward. The reader is referred to recent papers \cite{MTmultiD,KWo} for the specifics of the multi-depot VRP. For the open VRP (see references in recent papers \cite{LJ,LSHD}), it is enough to introduce a dummy depot with zero distances to this depot from all customers.

\item \emph{Tight capacity constraint.} For some instances of the VRP, the main difficulties lie in packing all goods into a bounded number of vehicles. Since the dynamic programming approach enumerates all possible subsets, the ``bin packing"/loading part of the VRP is resolved at the same time as the routing part. See arguments in \cite{CGS} for the benefits of integrating loading and routing. It is easy to see that the framework can incorporate  more complicated packing constraints, e.g. two-dimensional loading constraints \cite{LZZ}, by solving the corresponding packing sub-problems on each step of the calculations. 

\item \emph{Fixed items in a vehicle.} In some VRPs it is important to allocate customers to particular vehicles.
 We refer to these customers as \emph{fixed items} in a vehicle. This feature can easily be added to the dynamic programming recursions. Fixed items can be useful, e.g. with multiple visits of customers (see Section \ref{sec:solvable2TSP} in this paper). Another example is the so-called site-dependent VRP \cite{BaldachiVigo} - some customers can be served only by a specific type of vehicles (so-called \emph{docking} constraints \cite{CGRC}). An interesting case study with fixed customers was described in \cite{DGHL08}. An Austrian red cross considered introducing two tiers for a blood delivery service. Urgent delivery (with a higher price) is the delivery within one day, and standard delivery (at a lower price) is the delivery on the second day. Hospital customers for the current day are known, while the next day's customers are unknown and only the probabilities of requests can be evaluated. So, on each day a dispatcher knows undelivered requests from yesterday, requests which arrived today, and probable requests for tomorrow. The requests from yesterday are urgent and have to be delivered today - hence a fixed allocation of these customers to today's route (it is assumed here that the delivery is done by one vehicle). A sample of (probable) requests for tomorrow should be fixed for tomorrow's delivery. Two vehicles in the corresponding 2VRP are the vehicle for today's delivery and the vehicle for tomorrow's delivery. Today's requests are flexible and can be allocated to either vehicle, however they will be charged different prices.

\item \emph{Penalties for wrong day deliveries.} In the previous paragraph we mentioned charges/costs for deliveries in different vehicles/days. Another example is given in \cite{CRW15}, where penalties for wrong day deliveries were introduced. These can easily be incorporated into our model by changing  left/right lengths of intervals (assuming that costs of travelling and penalties are measured in the same monetary units).

\item \emph{Cumulative VRPs.} In the \emph{cumulative} VRP the objective is to minimise the sum of arrival times to all customers; this problem is also known as the \emph{latency} problem (see references in recent papers \cite{CSW,LW,RAP}). It is easy to see that the dynamic programming approach can be adapted for this type of objective function.

\item \emph{Weight or time dependent travel costs.}  There are some practical situations when the  travel costs depend on the load of the vehicle (see  \cite{ZTK,ZQZL} and references there) or on the time when the vehicle travels (\cite{ECT,TGJL}). It can easily be seen that the dynamic programming recursions above allow this type of calculation to be incorporated into the recursions.

\end{itemize}

The list above defines an exciting research project for the framework to be adopted and tested. In this paper we provide only limited however extremely impressive results of empirical tests of the framework on the instances of the balanced 2TSP discussed in Section~\ref{sec:solvable2TSP}.

\section
{\bf Empirical evaluation of the framework.}

\smallskip
The framework was implemented as the 2VRP\_Heuristic below. The main subroutine of the heuristic is the sliding subset procedure denoted as $H(s,l)$, where $s$ and $l$ are the two parameters of the framework. 
\begin{tabbing}
1111\=1111\=1111\=\kill
{\bf 2VRP\_Heuristic\{ }\\
\>{\bf repeat}  $m$ times:\\
\>\>GENERATE(next feasible solution);\\
\>\>\bf{repeat}\\
\>\>\>Apply  $H(s,l)$ for improving the feasible solution;\\
\>\>\>Apply BELPERM heuristic to improve individual tours \\
\>\>\>      and obtain in this way a new feasible solution;\\
\>\>\> Save the record;\\
\>\>{\bf until} (feasible solution cannot be improved);\\
\>{\bf return} record - the best solution found;\\
{\bf \}}
\end{tabbing}

We considered two modifications of the 2VRP\_Heuristic depending on the procedure used as GENERATE subroutine.
\smallskip

KSH-Heuristic uses KS heuristic (see Section~\ref{sec:solvable2TSP}) as the GENERATE subroutine. 
Number of repetitions $m$ is set to $n-1$, which is the number of possible start nodes in KS.  Next feasible solution on each step is defined by a new start point and by the corresponding Kalmanson sequence,  which is then used for permuting initial distance matrix. The permuted matrix is used  in recursions (\ref{eq:vi1})-(\ref{eq:vi3}) for obtaining the next feasible solution.
Computational experiments with the KS heuristic (see Section~\ref{sec:experiments1}, Table \ref{table:summaryHeuristics}) have shown a reasonably good performance of this heuristic. Therefore we use it here to get good initial feasible solutions. To understand the importance of good initial solutions for our framework, we have also considered random initial feasible solutions.

\smallskip
RndH heuristic uses  a simple generator of random feasible solutions as GENERATE subroutine. Number of repetitions  $m$ in RndH can be arbitrary big. In our computational experiments we set the number of repetitions to match the computational time of KSH heuristic with the same settings in $H(s,l)$ procedure.

\smallskip
We have tested these two heuristics on the two sets of instances: on the 2TSP benchmark instances from Bassetto \& Mason \cite{BassettoM11}, and on the set of randomly generated 2TSP instances with permuted Kalmanson matrices.
In our experiments, we used $H(3,1)$, $H(5,2)$ and $H(6,3)$ as subroutines. It means that on each step of local search we applied recursions (\ref{eq:dp0})-(\ref{eq:dpB}) for sub-problems with only 11, 15, or 17 nodes.
The results of the experiments are summarised in Tables~\ref{table:summary1}-\ref{table:summary3}. Computational times $t_m$ in the tables show the mean CPU time for the KSH heuristic (averaged over 20 instances). This time is used as stopping time for the RndH heuristic. Detailed results for all instances can be found in Appendix C.

\smallskip
Our computational experiments revealed two surprising results. Firstly, our approach has significantly improved the previously published results. We managed either to improve or match \emph{all} results obtained by computer algorithms and the majority of the results obtained with visualisation and manual intervention reported in \cite{BassettoM11}. 

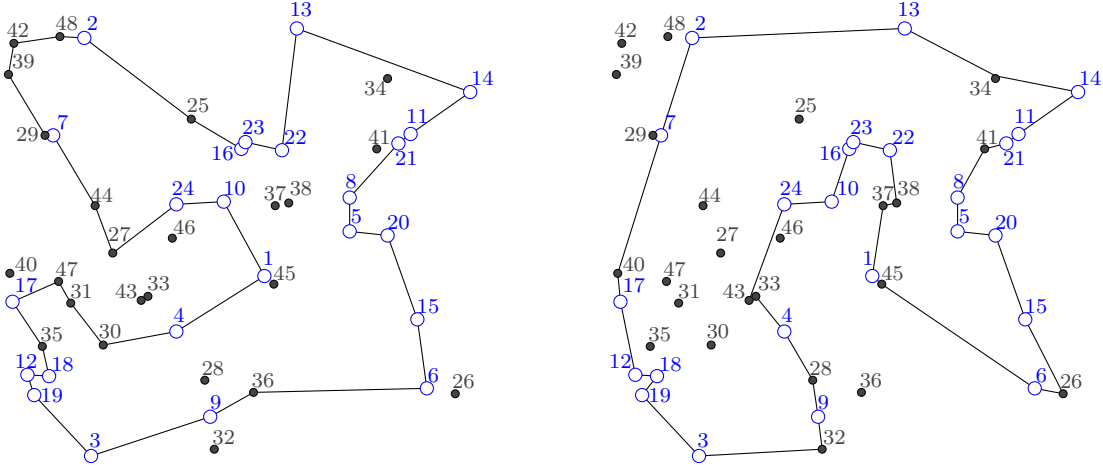
\begin{figure}
\unitlength=1cm
\definecolor{uuuuuu}{rgb}{0.26666666666666666,0.26666666666666666,0.26666666666666666}
\definecolor{qqqqff}{rgb}{0.,0.,1.}
\begin{tikzpicture}[line cap=round,line join=round,>=triangle 45,x=0.18cm,y=0.18cm]
\clip(0.,2.3) rectangle (39.911566481628,40.04844454925006);
\draw (21.3,18.1)-- (18.3,23.6)-- (14.8,23.4)-- (10.1,19.8)-- (8.8,23.3)-- (5.7,28.5)-- (5.1,28.5)-- (2.4,33.)-- (2.8,35.3)-- (6.2,35.8)-- (8.,35.7)-- (15.9,29.7)-- (19.6,27.5)-- (19.9,28.)-- (22.6,27.4)-- (23.7,36.4)-- (36.5,31.7)-- (32.1,28.6)-- (31.2,27.9)-- (27.6,23.9)-- (27.6,21.4)-- (30.4,21.1)-- (32.6,14.9)-- (33.3,9.8)-- (20.5,9.5)-- (17.3,7.7)-- (8.5,4.8)-- (4.3,9.3)-- (3.8,10.8)-- (5.4,10.7)-- (4.9,12.9)-- (2.7,16.2)-- (6.1,17.7)-- (7.,16.1)-- (9.4,13.)-- (14.8,14.);
\draw (21.3,18.1)-- (14.8,14.);
\begin{scriptsize}
\draw [color=qqqqff,fill=white] (21.3,18.1) circle (2.5pt);
\draw[color=qqqqff] (21.6,19.4) node {$1$};
\draw [color=qqqqff,fill=white] (8.,35.7) circle (2.5pt);
\draw[color=qqqqff] (8.6,36.8) node {$2$};
\draw [color=qqqqff,fill=white] (8.5,4.8) circle (2.5pt);
\draw[color=qqqqff] (8.8,6) node {3};
\draw [color=qqqqff,fill=white] (14.8,14.) circle (2.5pt);
\draw[color=qqqqff] (15.,15.2) node {$4$};
\draw [color=qqqqff,fill=white] (27.6,21.4) circle (2.5pt);
\draw[color=qqqqff] (28.1,22.5) node {$5$};
\draw [color=qqqqff,fill=white] (33.3,9.8) circle (2.5pt);
\draw[color=qqqqff] (33.8,10.8) node {$6$};
\draw [color=qqqqff,fill=white] (5.7,28.5) circle (2.5pt);
\draw[color=qqqqff] (6.4,29.5) node {$7$};
\draw [color=qqqqff,fill=white] (27.6,23.9) circle (2.5pt);
\draw[color=qqqqff] (27.7,25.1) node {$8$};
\draw [color=qqqqff,fill=white] (17.3,7.7) circle (2.5pt);
\draw[color=qqqqff] (17.7,8.8) node {$9$};
\draw [color=qqqqff,fill=white] (18.3,23.6) circle (2.5pt);
\draw[color=qqqqff] (19.1,24.6) node {$10$};
\draw [color=qqqqff,fill=white] (32.1,28.6) circle (2.5pt);
\draw[color=qqqqff] (32.6,29.9) node {$11$};
\draw [color=qqqqff,fill=white] (3.8,10.8) circle (2.5pt);
\draw[color=qqqqff] (3.7,11.8) node {$12$};
\draw [color=qqqqff,fill=white] (23.7,36.4) circle (2.5pt);
\draw[color=qqqqff] (24.,37.8) node {$13$};
\draw [color=qqqqff,fill=white] (36.5,31.7) circle (2.5pt);
\draw[color=qqqqff] (37.4,32.7) node {$14$};
\draw [color=qqqqff,fill=white] (32.6,14.9) circle (2.5pt);
\draw[color=qqqqff] (33.4,15.9) node {$15$};
\draw [color=qqqqff,fill=white] (19.6,27.5) circle (2.5pt);
\draw[color=qqqqff] (18.2,27.2) node {$16$};
\draw [color=qqqqff,fill=white] (2.7,16.2) circle (2.5pt);
\draw[color=qqqqff] (3.7,17.4) node {$17$};
\draw [color=qqqqff,fill=white] (5.4,10.7) circle (2.5pt);
\draw[color=qqqqff] (6.4,11.6) node {$18$};
\draw [color=qqqqff,fill=white] (4.3,9.3) circle (2.5pt);
\draw[color=qqqqff] (5.6,9.3) node {$19$};
\draw [color=qqqqff,fill=white] (30.4,21.1) circle (2.5pt);
\draw[color=qqqqff] (31.2,22.1) node {$20$};
\draw [color=qqqqff,fill=white] (31.2,27.9) circle (2.5pt);
\draw[color=qqqqff] (31.8,26.8) node {$21$};
\draw [color=qqqqff,fill=white] (22.6,27.4) circle (2.5pt);
\draw[color=qqqqff] (23.6,28.4) node {$22$};
\draw [color=qqqqff,fill=white] (19.9,28.) circle (2.5pt);
\draw[color=qqqqff] (20.5,29.1) node {$23$};
\draw [color=qqqqff,fill=white] (14.8,23.4) circle (2.5pt);
\draw[color=qqqqff] (15.3,24.6) node {$24$};
\draw [fill=uuuuuu] (15.9,29.7) circle (1.5pt);
\draw[color=uuuuuu] (16.4,30.7) node {$25$};
\draw [fill=uuuuuu] (35.4,9.4) circle (1.5pt);
\draw[color=uuuuuu] (36,10.3) node {$26$};
\draw [fill=uuuuuu] (10.1,19.8) circle (1.5pt);
\draw[color=uuuuuu] (10.6,21.2) node {$27$};
\draw [fill=uuuuuu] (16.9,10.4) circle (1.5pt);
\draw[color=uuuuuu] (17.5,11.3) node {$28$};
\draw [fill=uuuuuu] (5.1,28.5) circle (1.5pt);
\draw[color=uuuuuu] (3.8,28.6) node {$29$};
\draw [fill=uuuuuu] (9.4,13.) circle (1.5pt);
\draw[color=uuuuuu] (9.9,14) node {$30$};
\draw [fill=uuuuuu] (7.,16.1) circle (1.5pt);
\draw[color=uuuuuu] (7.8,17) node {$31$};
\draw [fill=uuuuuu] (17.6,5.3) circle (1.5pt);
\draw[color=uuuuuu] (18.3,6.2) node {$32$};
\draw [fill=uuuuuu] (12.7,16.6) circle (1.5pt);
\draw[color=uuuuuu] (13.5,17.6) node {$33$};
\draw [fill=uuuuuu] (30.4,32.7) circle (1.5pt);
\draw[color=uuuuuu] (29.5,32) node {$34$};
\draw [fill=uuuuuu] (4.9,12.9) circle (1.5pt);
\draw[color=uuuuuu] (5.6,13.8) node {$35$};
\draw [fill=uuuuuu] (20.5,9.5) circle (1.5pt);
\draw[color=uuuuuu] (21.2,10.4) node {$36$};
\draw [fill=uuuuuu] (22.1,23.3) circle (1.5pt);
\draw[color=uuuuuu] (22.1,24.2) node {$37$};
\draw [fill=uuuuuu] (23.1,23.5) circle (1.5pt);
\draw[color=uuuuuu] (24,24.6) node {$38$};
\draw [fill=uuuuuu] (2.4,33.) circle (1.5pt);
\draw[color=uuuuuu] (3.5,34) node {$39$};
\draw [fill=uuuuuu] (2.5,18.3) circle (1.5pt);
\draw[color=uuuuuu] (3.7,19.) node {$40$};
\draw [fill=uuuuuu] (29.6,27.5) circle (1.5pt);
\draw[color=uuuuuu] (29.8,28.4) node {$41$};
\draw [fill=uuuuuu] (2.8,35.3) circle (1.5pt);
\draw[color=uuuuuu] (3.1,36.3) node {$42$};
\draw [fill=uuuuuu] (12.2,16.3) circle (1.5pt);
\draw[color=uuuuuu] (11.1,17.) node {$43$};
\draw [fill=uuuuuu] (8.8,23.3) circle (1.5pt);
\draw[color=uuuuuu] (9.3,24.1) node {$44$};
\draw [fill=uuuuuu] (22.,17.5) circle (1.5pt);
\draw[color=uuuuuu] (22.8,18.3) node {$45$};
\draw [fill=uuuuuu] (14.5,20.9) circle (1.5pt);
\draw[color=uuuuuu] (15.4,21.7) node {$46$};
\draw [fill=uuuuuu] (6.1,17.7) circle (1.5pt);
\draw[color=uuuuuu] (6.7,18.6) node {$47$};
\draw [fill=uuuuuu] (6.2,35.8) circle (1.5pt);
\draw[color=uuuuuu] (6.7,36.8) node {$48$};
\end{scriptsize}
\end{tikzpicture}
\qquad
\begin{tikzpicture}[line cap=round,line join=round,>=triangle 45,x=0.18cm,y=0.18cm]
\clip(0.,2.3) rectangle (39.911566481628,39.94469485084606);
\draw (21.3,18.1)-- (22.,17.5)-- (33.3,9.8)-- (35.4,9.4)-- (32.6,14.9)-- (30.4,21.1)-- (27.6,21.4)-- (27.6,23.9)-- (29.6,27.5)-- (31.2,27.9)-- (32.1,28.6)-- (36.5,31.7)-- (30.4,32.9)-- (23.7,36.4)-- (8.,35.7)-- (5.7,28.5)-- (2.5,18.3)-- (2.7,16.2)-- (3.8,10.8)-- (5.4,10.7)-- (4.3,9.3)-- (8.5,4.8)-- (17.6,5.3)-- (17.3,7.7)-- (16.9,10.4)-- (14.8,14.)-- (12.7,16.6)-- (12.2,16.3)-- (14.8,23.4)-- (18.3,23.6)-- (19.6,27.5)-- (19.9,28.)-- (22.6,27.4)-- (23.1,23.5)-- (22.1,23.3);
\draw (22.1,23.3)-- (21.3,18.1);
\begin{scriptsize}
\draw [color=qqqqff,fill=white] (21.3,18.1) circle (2.5pt);
\draw[color=qqqqff] (21.0,19.4) node {$1$};
\draw [color=qqqqff,fill=white] (8.,35.7) circle (2.5pt);
\draw[color=qqqqff] (8.6,36.8) node {$2$};
\draw [color=qqqqff,fill=white] (8.5,4.8) circle (2.5pt);
\draw[color=qqqqff] (8.8,6) node {3};
\draw [color=qqqqff,fill=white] (14.8,14.) circle (2.5pt);
\draw[color=qqqqff] (15.,15.2) node {$4$};
\draw [color=qqqqff,fill=white] (27.6,21.4) circle (2.5pt);
\draw[color=qqqqff] (28.1,22.5) node {$5$};
\draw [color=qqqqff,fill=white] (33.3,9.8) circle (2.5pt);
\draw[color=qqqqff] (33.8,10.8) node {$6$};
\draw [color=qqqqff,fill=white] (5.7,28.5) circle (2.5pt);
\draw[color=qqqqff] (6.4,29.5) node {$7$};
\draw [color=qqqqff,fill=white] (27.6,23.9) circle (2.5pt);
\draw[color=qqqqff] (27.7,25.1) node {$8$};
\draw [color=qqqqff,fill=white] (17.3,7.7) circle (2.5pt);
\draw[color=qqqqff] (17.7,8.8) node {$9$};
\draw [color=qqqqff,fill=white] (18.3,23.6) circle (2.5pt);
\draw[color=qqqqff] (19.5,24.6) node {$10$};
\draw [color=qqqqff,fill=white] (32.1,28.6) circle (2.5pt);
\draw[color=qqqqff] (32.6,29.9) node {$11$};
\draw [color=qqqqff,fill=white] (3.8,10.8) circle (2.5pt);
\draw[color=qqqqff] (2.7,11.8) node {$12$};
\draw [color=qqqqff,fill=white] (23.7,36.4) circle (2.5pt);
\draw[color=qqqqff] (24.,37.8) node {$13$};
\draw [color=qqqqff,fill=white] (36.5,31.7) circle (2.5pt);
\draw[color=qqqqff] (37.4,32.7) node {$14$};
\draw [color=qqqqff,fill=white] (32.6,14.9) circle (2.5pt);
\draw[color=qqqqff] (33.4,15.9) node {$15$};
\draw [color=qqqqff,fill=white] (19.6,27.5) circle (2.5pt);
\draw[color=qqqqff] (18.2,27.2) node {$16$};
\draw [color=qqqqff,fill=white] (2.7,16.2) circle (2.5pt);
\draw[color=qqqqff] (3.7,17.4) node {$17$};
\draw [color=qqqqff,fill=white] (5.4,10.7) circle (2.5pt);
\draw[color=qqqqff] (6.4,11.6) node {$18$};
\draw [color=qqqqff,fill=white] (4.3,9.3) circle (2.5pt);
\draw[color=qqqqff] (5.6,9.3) node {$19$};
\draw [color=qqqqff,fill=white] (30.4,21.1) circle (2.5pt);
\draw[color=qqqqff] (31.2,22.1) node {$20$};
\draw [color=qqqqff,fill=white] (31.2,27.9) circle (2.5pt);
\draw[color=qqqqff] (31.8,26.8) node {$21$};
\draw [color=qqqqff,fill=white] (22.6,27.4) circle (2.5pt);
\draw[color=qqqqff] (23.6,28.4) node {$22$};
\draw [color=qqqqff,fill=white] (19.9,28.) circle (2.5pt);
\draw[color=qqqqff] (20.5,29.1) node {$23$};
\draw [color=qqqqff,fill=white] (14.8,23.4) circle (2.5pt);
\draw[color=qqqqff] (15.3,24.6) node {$24$};
\draw [fill=uuuuuu] (15.9,29.7) circle (1.5pt);
\draw[color=uuuuuu] (16.4,30.7) node {$25$};
\draw [fill=uuuuuu] (35.4,9.4) circle (1.5pt);
\draw[color=uuuuuu] (36,10.3) node {$26$};
\draw [fill=uuuuuu] (10.1,19.8) circle (1.5pt);
\draw[color=uuuuuu] (10.6,21.2) node {$27$};
\draw [fill=uuuuuu] (16.9,10.4) circle (1.5pt);
\draw[color=uuuuuu] (17.5,11.3) node {$28$};
\draw [fill=uuuuuu] (5.1,28.5) circle (1.5pt);
\draw[color=uuuuuu] (3.8,28.6) node {$29$};
\draw [fill=uuuuuu] (9.4,13.) circle (1.5pt);
\draw[color=uuuuuu] (9.9,14) node {$30$};
\draw [fill=uuuuuu] (7.,16.1) circle (1.5pt);
\draw[color=uuuuuu] (7.8,17) node {$31$};
\draw [fill=uuuuuu] (17.6,5.3) circle (1.5pt);
\draw[color=uuuuuu] (18.5,6.2) node {$32$};
\draw [fill=uuuuuu] (12.7,16.6) circle (1.5pt);
\draw[color=uuuuuu] (13.7,17.6) node {$33$};
\draw [fill=uuuuuu] (30.4,32.7) circle (1.5pt);
\draw[color=uuuuuu] (29.3,32) node {$34$};
\draw [fill=uuuuuu] (4.9,12.9) circle (1.5pt);
\draw[color=uuuuuu] (5.6,13.8) node {$35$};
\draw [fill=uuuuuu] (20.5,9.5) circle (1.5pt);
\draw[color=uuuuuu] (21.2,10.4) node {$36$};
\draw [fill=uuuuuu] (22.1,23.3) circle (1.5pt);
\draw[color=uuuuuu] (22.1,24.2) node {$37$};
\draw [fill=uuuuuu] (23.1,23.5) circle (1.5pt);
\draw[color=uuuuuu] (24,24.6) node {$38$};
\draw [fill=uuuuuu] (2.4,33.) circle (1.5pt);
\draw[color=uuuuuu] (3.5,34) node {$39$};
\draw [fill=uuuuuu] (2.5,18.3) circle (1.5pt);
\draw[color=uuuuuu] (3.7,19.) node {$40$};
\draw [fill=uuuuuu] (29.6,27.5) circle (1.5pt);
\draw[color=uuuuuu] (29.8,28.4) node {$41$};
\draw [fill=uuuuuu] (2.8,35.3) circle (1.5pt);
\draw[color=uuuuuu] (3.1,36.3) node {$42$};
\draw [fill=uuuuuu] (12.2,16.3) circle (1.5pt);
\draw[color=uuuuuu] (11.1,17.) node {$43$};
\draw [fill=uuuuuu] (8.8,23.3) circle (1.5pt);
\draw[color=uuuuuu] (9.3,24.1) node {$44$};
\draw [fill=uuuuuu] (22.,17.5) circle (1.5pt);
\draw[color=uuuuuu] (22.8,18.3) node {$45$};
\draw [fill=uuuuuu] (14.5,20.9) circle (1.5pt);
\draw[color=uuuuuu] (15.4,21.7) node {$46$};
\draw [fill=uuuuuu] (6.1,17.7) circle (1.5pt);
\draw[color=uuuuuu] (6.7,18.6) node {$47$};
\draw [fill=uuuuuu] (6.2,35.8) circle (1.5pt);
\draw[color=uuuuuu] (6.7,36.8) node {$48$};
\end{scriptsize}
\end{tikzpicture}
\caption{Solution to instance 55 with the total length of two tours of 33507; Nodes $\{1,2,\ldots,24\}$ appear in the both tours. Solution is found by RndH with H(3,1) settings after generating 63 initial random solution.}
\label{fig:inst55}
\end{figure}

Within the described settings of our experiments, only the solution found for instance 55 from \cite{BassettoM11} was worse than the previously known solution. We did manage to match the best known solution to this instance by changing the settings and increasing the number of repetitions in RndH to 63 random starts. The solution is depicted in Fig.~\ref{fig:inst55}. It is difficult to say what made this instance the most difficult one for our algorithms.

\begin{table}
\begin{center}
\begin{tabular}{||c||l||c|c||c|c||c|c||}
\hline 
  \multicolumn{2}{||c||}{} &
  \multicolumn{2}{c||}{$H(3,1)$, $t_m=18$s  } & \multicolumn{2}{c||}{$H(5,2)$, $t_m=134$s } & \multicolumn{2}{c||}{$H(6,3)$, $t_m=334$s  } \\ 
\cline{3-8} 
\multicolumn{2}{||c||}{}  & PC &PC+m & PC & PC+m & PC & PC+m \\ 
\hline
KSH&Mean \% & $-$2.29  & $-$0.58  & $-$2.73  & $-$1.03  & $-$2.63  & $-$0.93  \\ 
& Best \%& $-$6.77  & $-$1.99  & $-$7.31  & $-$2.77  & $-$7.31  & $-$2.77  \\ 
&Worst \%& $-$0.41  & +1.25  & $-$0.52  & +0.00  & $-$0.52  & +0.00  \\ 
&Improved \#& 20/20 & 15/20 & 20/20 & 20/20 & 20/20 & 20/20 \\ 
\hline 
RndH&Mean \% & $-$2.05  & $-$0.33  & $-$2.55  & $-$0.84  & $-$2.63  & $-$0.93  \\ 
& Best \%& $-$6.98  & $-$1.83  & $-$6.92  & $-$2.77  & $-$7.19  & $-$2.77  \\ 
&Worst \%& $-$0.39  & +0.99  & $-$0.32  & +0.78  & $-$0.39  & +0.15  \\ 
&Improved \#& 20/20 & 12/20 & 20/20 & 18/20 & 20/20 & 19/20 \\ 
\hline
\end{tabular} 
\end{center}
\caption{Summary of results for instances with $8$ (out of $48$) nodes  visited in two periods.   }
\label{table:summary1}
\end{table}


\begin{table}
\begin{center}
\begin{tabular}{||c||l||c|c||c|c||c|c||}
\hline 
  \multicolumn{2}{||c||}{} &
  \multicolumn{2}{c||}{$H(3,1)$, $t_m=18$s  } & \multicolumn{2}{c||}{$H(5,2)$, $t_m=133$s } & \multicolumn{2}{c||}{$H(6,3)$, $t_m=395$s  } \\ 
\cline{3-8} 
\multicolumn{2}{||c||}{}  & PC &PC+m & PC & PC+m & PC & PC+m \\ 
\hline
KSH&Mean \% & $-$1.71  & $-$0.68  & $-$2.14  & $-$1.11  & $-$2.17  & $-$1.14  \\ 
& Best \%& $-$3.68  & $-$1.98  & $-$3.93  & $-$2.59  & $-$3.93  & $-$2.59  \\ 
&Worst \%& +0.04  & +1.44  & $-$0.52  & +0.47  & $-$0.54  & +0.51  \\ 
&Improved \#& 19/20 & 15/20 & 20/20 & 18/20 & 20/20 & 19/20 \\ 
\hline 
RndH&Mean \% & $-$1.33  & $-$0.29  & $-$1.97  & $-$0.94  & $-$2.03  & $-$1.00  \\ 
& Best \%& $-$3.19  & $-$1.98  & $-$3.93  & $-$2.59  & $-$3.86  & $-$2.52  \\ 
&Worst \%& +0.69  & +2.86  & $-$0.45  & +0.93  & $-$0.52  & +0.39  \\ 
&Improved \#& 18/20 & 12/20 & 20/20 & 18/20 & 20/20 & 18/20 \\ 
\hline
\end{tabular} 
\end{center}
\caption{Summary of results for instances with $16$ (out of $48$) nodes  visited in two periods. }
\label{table:summary2}
\end{table}

\begin{table}
\begin{center}
\begin{tabular}{||c||l||c|c||c|c||c|c||}
\hline 
  \multicolumn{2}{||c||}{} &
  \multicolumn{2}{c||}{$H(3,1)$, $t_m=14$s  } & \multicolumn{2}{c||}{$H(5,2)$, $t_m=119$s } & \multicolumn{2}{c||}{$H(6,3)$, $t_m=255$s  } \\ 
\cline{3-8} 
\multicolumn{2}{||c||}{}  & PC &PC+m & PC & PC+m & PC & PC+m \\ 
\hline
KSH&Mean \% & $-$1.12  & $-$0.50  & $-$1.31  & $-$0.69 & $-$1.32 & $-$0.70 \\ 
&  Best \%& $-$3.70 & $-$1.52 & $-$4.01 & $-$1.71 & $-$4.01 & $-$1.60 \\
&Worst \%& +0.00 & +0.40 & $-$0.17 & +0.24 & $-$0.17 & +0.33 \\
&Improved \#& 20/20 & 15/20 & 20/20 & 17/20  & 20/20 & 17/20 \\ 
\hline 
RndH&Mean \% & $-$0.69 & $-$0.07 & $-$1.25 & $-$0.63 & $-$1.16 & $-$0.50 \\ 
&  Best \%& $-$3.68 & $-$1.42 & $-$3.88 & $-$2.17 & $-$3.98 & $-$1.97 \\  
&Worst \%& +0.84 & +1.07 & $-$0.16 & +0.17 & +0.07 & +1.04 \\  
&Improved \#& 17/20 & 11/20 & 20/20 & 17/20 & 19/20 & 13/20 \\ 
\hline
\end{tabular} 
\end{center}
\caption{Summary of results for instances with $24$ (out of $48$) nodes  visited in two periods.}
\label{table:summary3}
\end{table}

\smallskip
The second surprising result in our experiments is the fact that the quality of an initial solution is not crucial for the quality of the final results. As one can see, the quality of the results for the RndH heuristic, where initial solutions are randomly generated, are almost as good as the results for the KSH.

\smallskip
{\bf Experiments with random Kalmanson matrices.} To further test the computational efficiency of our framework, we generated instances with random permuted strong Kalmanson matrices. For this, we used an approach discussed in Section~\ref{sec:solvable2TSP}. We generated all random numbers used, i.e.\ items in the first row and column  as well as $\alpha_{ij}$ and $\beta_i$, from the interval $[0.1,1.1]$ by using $C$-build-in function $rand()$. Two sets of instances have been generated. Each set contains 60 instances. The size of instances in the first set is $50/30$, meaning that we have $50$ points with $30$ points to be visited twice. The size of instances in the second set is $100/30$. Since the distance matrices in the instances are permuted Kalmanson matrices, there is no sense to test the KSH heuristic (the heuristic would find optimal solutions for all instances). So, we tested the RndH heuristic with various settings for the framework. To manage the time spent on computational experiments we excluded $H(6,3)$ settings, and compared the settings $H(2,1)$, $H(3,1)$, and $H(5,2)$. The number of iterations was set to $100$, and we recorded performance for $1$, $5$, $10$, $20$, $30$, $40$, $50$, $60$, $70$, $80$, $90$ and $100$ iterations. 

\begin{figure}
\unitlength=1cm
{
\begin{picture}(10.5,8)
\includegraphics[scale=1]{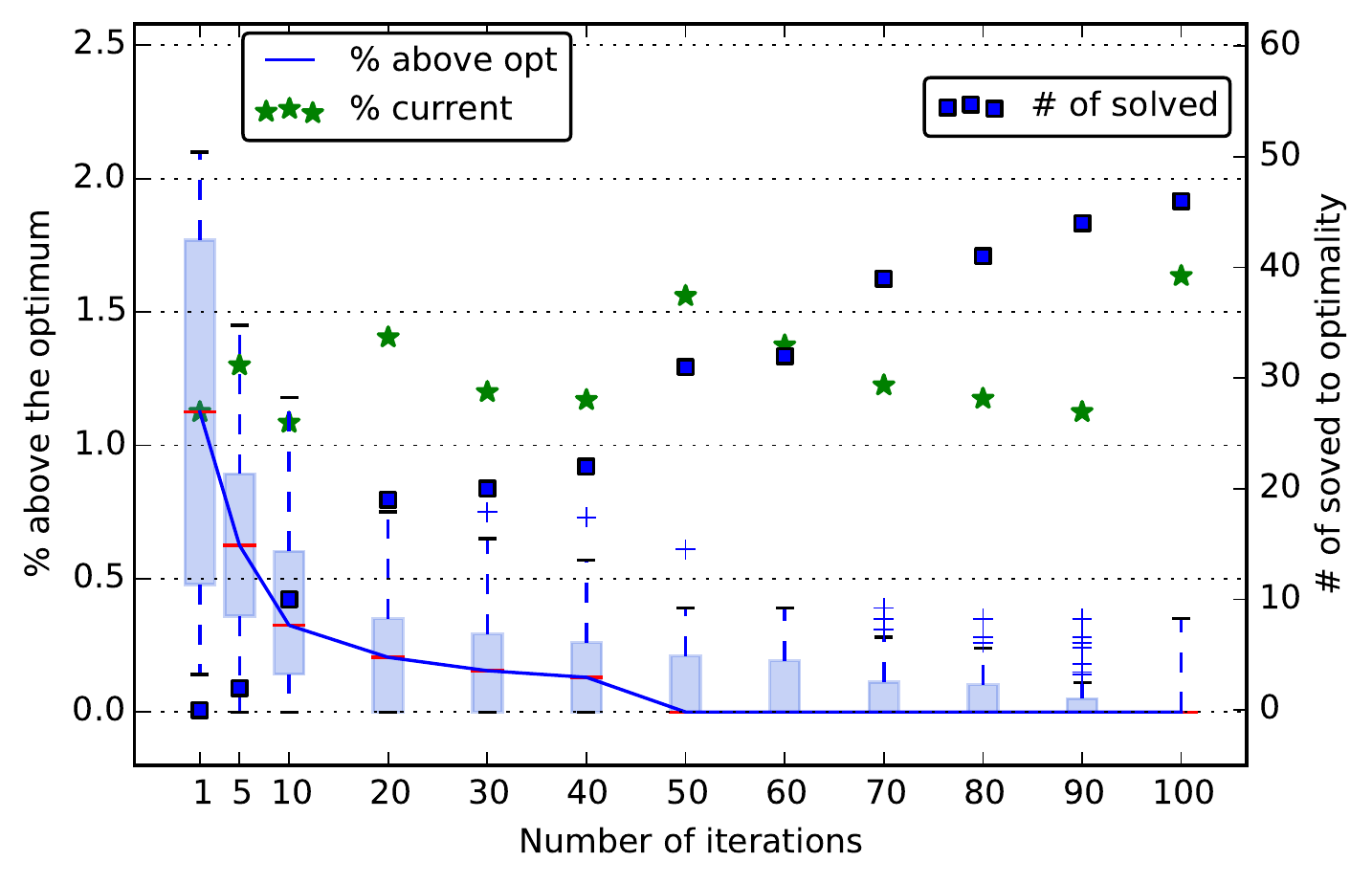}
\end{picture} 
}
\caption{$H(3,1)$ settings: Summary of results for the $50/30$-set of $60$ Kalmanson instances.}
\label{fig:in50_31}
\end{figure}

Graphs in Fig.~\ref{fig:in50_31} show the results of computational experiments for the set $50/30$ and $H(3,1)$ setting for the framework. For the chosen number of iterations we  recorded percentage above the optimum for each of $60$ instances. The graph in the figure shows the box-\&-whiskers plot for this set of results. As we can see, the median after the first iteration is just above $1\%$, while after $50$ iterations it is already $0$. We have also recorded the number of instances solved to optimality: these are  $0,2,10,19,20,22,31,32,39,41,44,46$. So after $100$ iterations, $46$ instances have been solved to optimality. 
To get an idea about accuracy of the RndH procedure on each iteration, we recorded the mean for solutions obtained on each iteration (``$*$" symbols on the graph). For instance, solutions for randomly generated initial allocations on the last step were the worst according to the graph. This indicates that it is important to have many random starts to eventually get very good solutions.

The results are even more impressive for the $H(5,2)$ settings of the framework (see charts in Appendix D): after $60$ iterations \emph{all} $60$ instances have been solved to optimality! 

The size of instances in the $50/30$-set of test problems is similar to the size of instances from \cite{BassettoM11}: each vehicles visits $40$ customers. On both sets of instance the framework has shown impressively good performance. Will an increase in the size of instance influence the performance of the framework?

In the second set of random instances the size of test problems is slightly bigger: each vehicle visits now $(100+30)/2=65$ customers. 

\begin{figure}
\unitlength=1cm
{
\begin{picture}(10.5,8)
\includegraphics[scale=1]{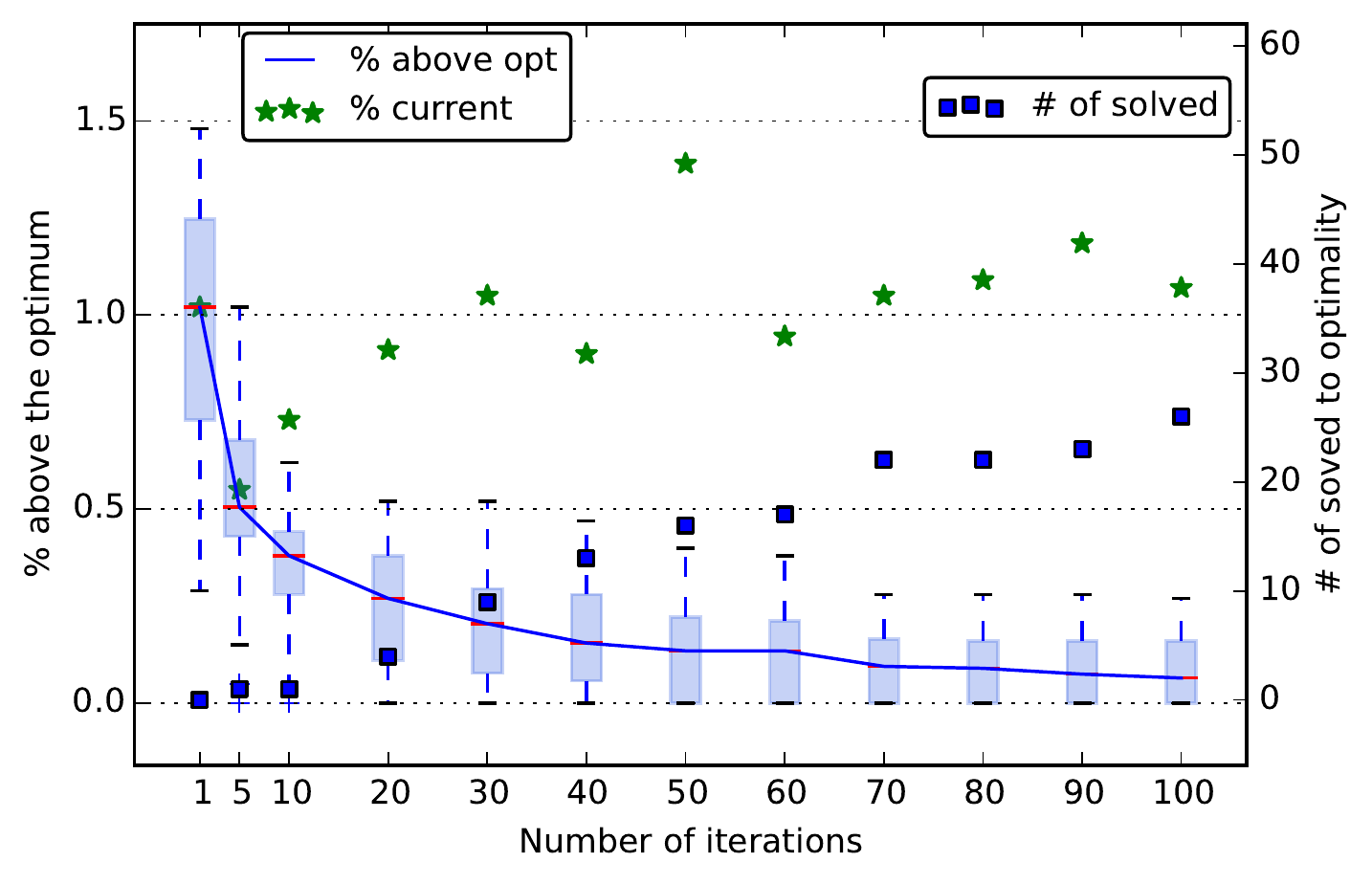}
\end{picture} 
}
\caption{RndH heuristic with $H(5,2)$: Summary for the $100/30$ set of $60$ Kalmanson instances.}
\label{fig:in100_52}
\end{figure}

Results for the $H(5,2)$ settings are shown in Fig.~\ref{fig:in100_52}. After 100 iterations, 26 instances have been solved to optimality.  
The number of instances solved to optimality for the chosen number of iterations are 0,	1,	1,	4,	9,	13,	16,	17,	22, 22,	23, and 26, correspondingly. 
Appendix E shows summary of results for three settings of the framework in the RndH algorithm. For the $H(2,1)$ setting, none of the instances was solved to optimality. For the $H(3,1)$ settings, after $50$ iterations only $2$ instances have been solved to optimality; after $100$ iterations only $5$ problems have been solved to optimality. 
As one can guess, the improvement in the accuracy achieved with $H(5,2)$ settings compared to $H(3,1)$ settings came at a price: the price is the CPU time. 

\begin{figure}
\unitlength=1cm
{
\begin{picture}(10.5,8)
\includegraphics[scale=1]{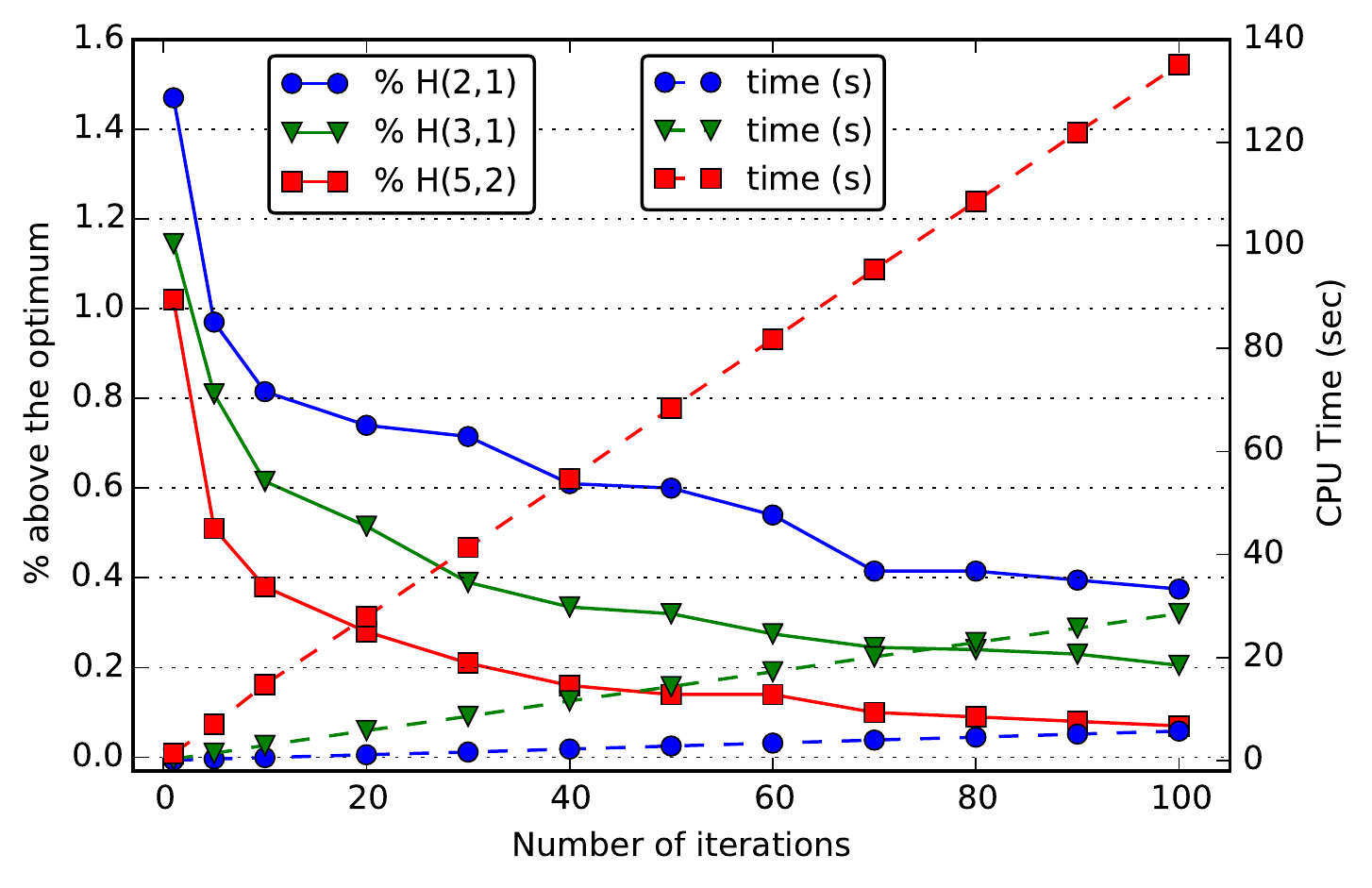}
\end{picture} 
}
\caption{Comparative analysis of settings $H(2,1)$, $H(3,1)$, and $H(5,2)$ on the 100/30 set.}
\label{fig:in100_time}
\end{figure}

\begin{figure}
\unitlength=1cm
{
\begin{picture}(10.5,8)
\includegraphics[scale=1]{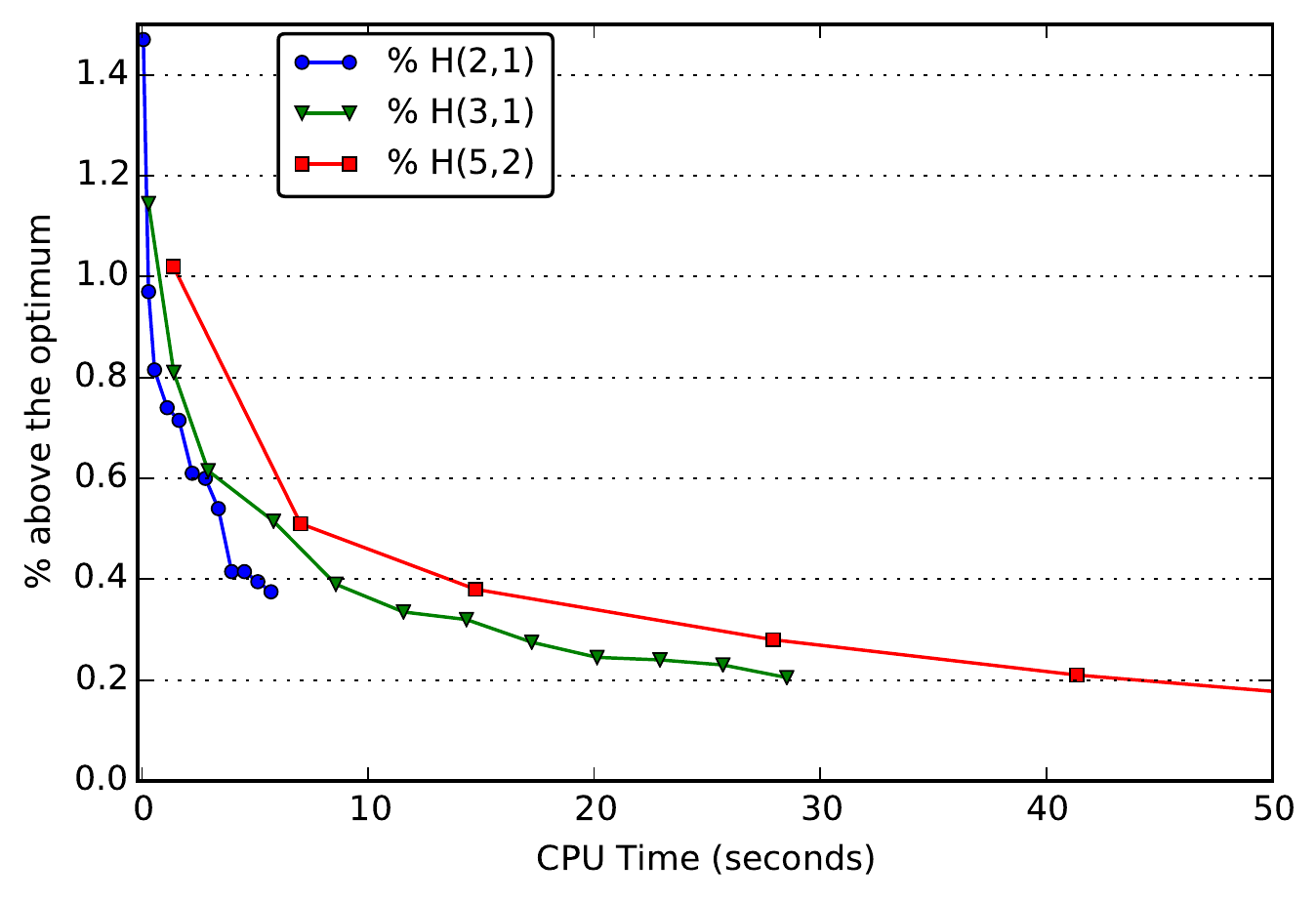}
\end{picture} 
}
\caption{Comparative analysis of accuracy vs CPU time for $H(2,1)$, $H(3,1)$, and $H(5,2)$ settings}
\label{fig:in100_Times}
\end{figure}

Fig.~\ref{fig:in100_time} compares computational time and accuracy for three chosen settings of the framework for 100/30 set. The main axes on the plot shows the mean of percentage above the optimum, the secondary axis shows the CPU time. Fig.~\ref{fig:in100_Times}, where the percentage above the optimum plotted against the time spent, clearly illustrates that none of the settings dominates the others.

\section{Summary} In this paper we have presented a framework for small but rich VRPs. The framework can be used for a variety of the VRP settings. The computational scheme which implements the framework permits an easy scaling,  hence varying computational time and the accuracy of the solutions found. We have tested our approach on the balanced 2TSP - a special case of the VRP with only two vehicles (2VRP). We have described a new polynomially solvable case of the 2TSP with the Kalmanson distance matrices.

Computational experiments on the previously published test instances for the balanced 2TSP have shown an impressive performance of the framework. We managed to improve solutions for 55 out of 60 benchmark instances. For the remaining 5 instances we found solutions with the same value of the objective function. In the absence of alternative sets of test instances, we decided to use easy solvable cases of the 2TSP for generating new sets of test instances. For the set of new instances of about the same size of problems as in our previous experiments, the performance of the framework was also very good: all generated instances had been solved to optimality. However, increasing the size of the problems from 80 to 130 revealed, as was expected, an imperfection of the framework: only 26 out of 60 instances had been solved to optimality. Thus, from methodological point of view the use of easy solvable cases for generating test instances appeared to be productive.

The framework described in the paper can be easily extended to the case with more than 2 vehicles, e.g. as a simple heuristic which considers all possible pairs of vehicles  and the corresponding 2VRPs. We have successfully used this approach in real life applications (to be reported in our next publications) as well as in international contents in logistic optimisations. Our computational results where we used the framework, have got the third prize in the VEROLOG competition in 2015 (Vienna, VEROLOG-2015).

\section*{Acknowledgements}
\nopagebreak
Financial support by Austrian Science Fund (FWF): W1230, Doctoral Program
''Discrete Mathematics'' and by the Center for Discrete Mathematics and Its 
Applications, University of Warwick, is acknowledged. The authors thank 
Tatiana Bassetto and Francesco Mason for providing the benchmark problems as
well as Douglas Miranda and Alex Taylor for useful comments on an early 
version of the paper.

\pagebreak
\appendix
\section{Appendix: $O(n^2)$ space recursions for the balanced 2TSP}

The algorithm starts with defining boundary conditions for the recursions.

\begin{equation}
\begin{aligned}
MV(i,m)=
\begin{cases}
c_{i1}+c_{1n},\ \ \mbox{if}\  m=k,\\
\infty,\ \mbox{otherwise;}
\end{cases}\\
MW(i,m)=\begin{cases}
c_{n1}+c_{1i},\ \ \mbox{if}\  m=k,\\
\infty,\  \mbox{otherwise;}
\end{cases}\\
i=1,2,\ldots,n-1.\\
D(n,m)=
\begin{cases}
c_{n1}+c_{1n},\ \ \mbox{if}\  n\in S,\\
\infty,\  \mbox{otherwise.}
\end{cases}
\end{aligned}
\end{equation}

Then the values below are calculated for $j=n-1,n-2,\ldots,2$ :
\begin{equation}
\begin{aligned}
MV(i,m)=\begin{cases}
min\begin{cases}
c_{i,j+1}+MW(j,m)\\
c_{j+1,j}+MV(i,m+1),\ \ \ \ \ \mbox{if}\ j+1\notin S;
\end{cases}\\
c_{i,j+1}+c_{j+1,j}+D(j+1,m+1),\ \mbox{if } j+1\in S;\\
\hspace{3cm} i=1,2,\ldots,j-1.
\end{cases}
\end{aligned}
\end{equation}

\begin{equation}
\begin{aligned}
MW(i,m)=\begin{cases}
min\begin{cases}
c_{j,j+1}+MW(i,m)\\
c_{j+1,i}+MV(j,m+1),\ \ \ \ \ \mbox{if}\ j+1\notin S;
\end{cases}\\
c_{j,j+1}+c_{j+1,j}+D(j+1,m+1),\ \mbox{if } j+1\in S;\\
\hspace{3cm} i=1,2,\ldots,j-1.
\end{cases}
\end{aligned}
\end{equation}

\begin{equation}
\begin{aligned}
D(j,m)=\begin{cases}
min\begin{cases}
c_{j,j+1}+MW(j,m)\\
c_{j+1,j}+MV(j,m+1),\ \ \ \ \ \mbox{if}\ j+1\notin S;
\end{cases}\\
c_{j,j+1}+c_{j+1,j}+D(j+1,m+1),\ \mbox{if } j+1\in S.
\end{cases}
\end{aligned}
\end{equation}

And eventually the optimal length is calculated as 
\begin{equation}
\begin{aligned}
F_{opt}=\begin{cases}
min\begin{cases}
c_{12}+MW(1,1)\\
c_{21}+MV(1,2),\ \ \ \ \ \mbox{if}\ j+1\notin S;
\end{cases}\\
c_{12}+c_{21}+D(j22),\ \mbox{if } 2\in S.
\end{cases}
\end{aligned}
\end{equation}

\section {Appendix: Results of tests for KS-Heuristic (PC vs KS) }\nopagebreak[4]
\begin{small}
\begin{center}
\begin{tabular}{||c|c|c||c||}
\multicolumn{4}{l}{{Tests for the set with {\bf 8} nodes visited twice}}\\
\hline 
  \multirow{2}{*}& \multicolumn{2}{c||}{Solutions from \cite{BassettoM11}}& KS   \\ 
\cline{2-4}
Instance & PC & PC+manual &length  \\ 
\hline 
$I_{21}$	& 25217	& 24937	&  \underline{25162}\\
$I_{22}$	& \underline{26996}	& 26549		& 27564\\
$I_{23}$	& \underline{ 26476}	& 26192		& 26597\\
$I_{24}$	& 26802	& 26038	& \underline{26210}\\
$I_{25}$	& 27728	& 27408		& \underline{27488}\\
$I_{26}$	& \underline{24348}	& 24268		& 24473\\
$I_{27}$	& \underline{27335}	& 26857		& 27523\\
$I_{28}$	& \underline{24679}	& 24232		& 24712\\
$I_{29}$	& 26890	& 26466		& \underline{26690}\\
$I_{30}$	& 24978	& 24200		& \underline{24967}\\
$I_{31}$	& 26266	& 26130		& \underline{26138}\\
$I_{32}$	& 26360	& 26054		& \underline{ 26287}\\
$I_{33}$	& \underline{26418}	& 26418		& 26947\\
$I_{34}$	& 28733	& 27074		& \underline{27166}\\
$I_{35}$	& 25043	& 24587		& \underline{24706}\\
$I_{36}$	& 27103	& 26790		& \underline{26603}\\
$I_{37}$	& 25662	& 25123		& \underline{24893}\\
$I_{38}$	& 26459	& 25709		& \underline{26005}\\
$I_{39}$	& 27209	& 26994		& \underline{27099}\\
$I_{40}$	& \underline{25416}	& 24964		& \underline{25416}\\
\hline 
\end{tabular} 
\qquad
\begin{tabular}{||c|c|c||c||}
\multicolumn{4}{l}{{Tests for the set with {\bf 16} nodes visited twice}}\\
\hline  
  \multirow{2}{*}& \multicolumn{2}{c||}{Solutions from \cite{BassettoM11}}  & KS   \\ 
\cline{2-4} 
Instance & PC & PC+manual &length   \\ 
\hline 
$I_{1}$	& 33804	& 32556		& \underline{32888}\\
$I_{2}$	& \underline{30929}	& 30929		& 31898\\
$I_{3}$	& 30596	& 30382		&  \underline{30381}\\
$I_{4}$	& \underline{28563}	& 28441		& 28915\\
$I_{5}$	& \underline{27323}	& 27323		& 27594\\
$I_{6}$	& 33065	& 32546	&  \underline{ 32781}\\
$I_{7}$	& 32854	& 31861	&  \underline{ 32253}\\
$I_{8}$	& \underline{30850}	& 30571		& 30887\\
$I_{9}$	& \underline{34709}	& 34024		& 34864\\
$I_{10}$	& 31451	& 31201&  \underline{ 31265}\\
$I_{61}$	& 27158	& 26934		& \underline{26904}\\
$I_{62}$	& 27774	& 27619		& \underline{27735}\\
$I_{63}$	& 25308	& 24960		& \underline{25122}\\
$I_{64}$	& 27875	& 27285		&  \underline{ 27793}\\
$I_{65}$	& 27060	& 26888		& \underline{26885}\\
$I_{66}$	& 27677	& 27624		&  \underline{27339}\\
$I_{67}$	& 30268	& 30203		&  \underline{30028}\\
$I_{68}$	& 28033	& 27923		& \underline{27958}\\
$I_{69}$	& 27958	& 27638	&  \underline{27466}\\
$I_{70}$	& 28483	& 28427		& \underline{28218}\\
\hline 
\end{tabular} 
\end{center}
\begin{center}
\begin{tabular}{||c|c|c||c||}
\multicolumn{4}{l}{{Tests for the set with {\bf 24} nodes visited twice}}\\
\hline 
  \multirow{2}{*}& \multicolumn{2}{c||}{Solutions from \cite{BassettoM11}}& KS   \\ 
\cline{2-4}
Instance & PC & PC+manual &length  \\ 
\hline 
$I_{41}$	&  \underline{30253}	& 30147		&  \underline{30253}\\
$I_{42}$	& 33008	& 32020	&  \underline{32413}\\
$I_{43}$	& 31500	& 31500		&  \underline{31426}\\
$I_{44}$	& \underline{30313}	& 30170		& 30328\\
$I_{45}$	& 27986	& 27857		& \underline{27969}\\
$I_{46}$	&  \underline{30073}	& 30013		& 30368\\
$I_{47}$	& 32106	& 32106	&  \underline{31790}\\
$I_{48}$	& 31004	& 30942		& \underline{30933}\\
$I_{49}$	& \underline{33663}	& 33185		& 33818\\
$I_{50}$	& 31266	& 31266	&  \underline{31212}\\
$I_{51}$	& 33722	& 33627		&  \underline{33635}\\
$I_{52}$	& 32353	& 32280		&  \underline{32305}\\
$I_{53}$	& 33287	& 33200		&  \underline{32727}\\
$I_{54}$	& 31973	& 31600		&  \underline{31581}\\
$I_{55}$	& 33837	& 33507		&  \underline{33832}\\
$I_{56}$	&  \underline{29696}		& 30109	& 29762\\
$I_{57}$	& 31954	& 31640	& \underline{31827}\\
$I_{58}$	& 30705	& 30246		&  \underline{30372}\\
$I_{59}$	& 31549	& 31549	&  \underline{31223}\\
$I_{60}$	&  \underline{32384}		& 32839	& 32749\\
\hline 
\end{tabular} 
\end{center}
\end{small}

\section {Appendix: Results of tests for KSH and RndH heuristics } \nopagebreak[4]
\begin{small}

\begin{center}
\begin{tabular}{||c|c|c||c|c|c||c|c|c||c|c|c||}
\multicolumn{12}{l}{Tests for the set with 8 nodes visited twice}\\
\hline 
  \multirow{2}{*}& \multicolumn{2}{c||}{Results from \cite{BassettoM11}}& \multicolumn{3}{c||}{$H(3,1)$, $t_{mean}=18$s } & \multicolumn{3}{c||}{$H(5,2)$, $t_{mean}=134$s }  & \multicolumn{3}{c||}{$H(6,3)$, $t_{mean}=334$s } \\ 
\cline{2-3} \cline{4-6} \cline{7-9}  \cline{10-12} 
Ins & PC & PC+m & $t_{KSH}$ & KSH & RndH &$t_{KSH}$ & KSH & RndH & $t_{KSH}$& KSH& RndH \\ 
\hline 
$I_{21}$	& {25217}	& {24937}	& 18	& {24565}	& {24554}	& 132	& {24523}	& {24517}	& 320	& {24517}	& \underline{24493}\\
$I_{22}$	& {26996}	& {26549}	& 22	& {26456}	& {26430}	& 155	& \underline{26347}	& {26361}	& 382	& {26361}	& {26361}\\
$I_{23}$	& {26476}	& {26192}	& 19	& {26215}	& {26300}	& 137	& \underline{26050}	& {26180}	& 344	& {26050}	& {26050}\\
$I_{24}$	& {26802}	& {26038}	& 16	& {26117}	& {26116}	& 118	& {26032}	& \underline{26002}	& 304	& {26002}	& {26002}\\
$I_{25}$	& {27728}	& {27408}	& 21	& {26863}	& {27196}	& 134	& {26850}	& \underline{26790}	& 348	& {26843}	& {26914}\\
$I_{26}$	& {24348}	& {24268}	& 18	& {23882}	& {23982}	& 140	& \underline{23836}	& {23836}	& 296	& {23836}	& {23882}\\
$I_{27}$	& {27335}	& {26857}	& 16	& {26934}	& {26909}	& 127	& {26566}	& \underline{26417}	& 349	& {26712}	& {26663}\\
$I_{28}$	& {24679}	& {24232}	& 21	& {24028}	& {24010}	& 140	& \underline{23936}	& {23961}	& 381	& {23936}	& {23936}\\
$I_{29}$	& {26890}	& \underline{26466}	& 18	& \underline{26466}	& {26569}	& 123	& {26466}	& {26672}	& 322	& {26466}	& {26506}\\
$I_{30}$	& {24978}	& {24200}	& 20	& {23979}	& {24177}	& 145	& \underline{23915}	& {23915}	& 388	& {23915}	& {23915}\\
$I_{31}$	& {26266}	& \underline{26130}	& 16	&\underline {26130}	& {26158}	& 130	& {26130}	& {26130}	& 300	& {26130}	& {26130}\\
$I_{32}$	& {26360}	& {26054}	& 16	& {26072}	& {26117}	& 136	& \underline{26032}	& {26058}	& 346	& {26032}	& {26032}\\
$I_{33}$	& {26418}	& {26418}	& 20	& {26309}	& {26315}	& 142	&\underline{26216}	& {26333}	& 342	& {26251}	& {26315}\\
$I_{34}$	& {28733}	& {27074}	& 21	& {26788}	& {26728}	& 153	& \underline{26633}	& {26743}	& 354	& {26633}	& {26666}\\
$I_{35}$	& {25043}	& {24587}	& 15	& \underline{24516}	& {24576}	& 112	& {24516}	& {24516}	& 295	& {24516}	& {24516}\\
$I_{36}$	& {27103}	& {26790}	& 19	& {26391}	& {26657}	& 136	& {26079}	& {26368}	& 353	& {26285}	& \underline{26078}\\
$I_{37}$	& {25662}	& {25123}	& 16	& \underline{24893}	& {25136}	& 104	& {24893}	& {25072}	& 262	& {24893}	& {25070}\\
$I_{38}$	& {26459}	& {25709}	& 18	& \underline{25587}	& {25963}	& 127	& {25587}	& {25624}	& 314	& {25587}	& {25587}\\
$I_{39}$	& {27209}	& {26994}	& 19	& {26483}	& {26500}	& 156	& \underline{26246}	& {26246}	& 380	& {26246}	& {26246}\\
$I_{40}$	& {25416}	& {24964}	& 16	& {25276}	& {24861}	& 128	&\underline {24786}	& {24861}	& 301	& {24951}	& {24786}\\
\hline 
\end{tabular} 
\end{center}
\begin{center}
\begin{tabular}{||c|c|c||c|c|c||c|c|c||c|c|c||}
\multicolumn{12}{l}{Tests for the set with 16 nodes visited twice}\\
\hline 
  \multirow{2}{*}& \multicolumn{2}{c||}{Results from \cite{BassettoM11}}& \multicolumn{3}{c||}{$H(3,1)$, $t_{mean}=18$s} & \multicolumn{3}{c||}{$H(5,2)$, $t_{mean}=133$s} & \multicolumn{3}{c||}{$H(6,3)$, $t_{mean}=395$s } \\ 
\cline{2-3} \cline{4-6} \cline{7-9}  \cline{10-12} 
Ins & PC & PC+m &$t_{KSH}$ & KSH & RndH &$t_{KSH}$ & KSH&RndH &$t_{KSH}$&KSH&RndH \\ 
\hline 
$I_{1}$	& {33804}	& {32556}	& 17	& {32561}	& {32755}	& 116	& {32561}	& \underline{32555}	& 313	& {32555}	& {32555}\\
$I_{2}$	& {30929}	& {30929}	& 20	& {30940}	& {30811}	& 147	& {30768}	& {30768}	& 435	& \underline{30744}	& {30768}\\
$I_{3}$	& {30596}	& {30382}	& 16	& {30260}	& {30436}	& 116	& \underline{30086}	& {30104}	& 363	& {30086}	& {30098}\\
$I_{4}$	& {28563}	& {28441}	& 20	& {28353}	& {28274}	& 156	& \underline{28223}	& {28297}	& 398	& {28260}	& {28260}\\
$I_{5}$	& {27323}	& {27323}	& 18	& {27187}	& {27209}	& 153	& \underline{27176}	& {27176}	& 433	& {27176}	& {27176}\\
$I_{6}$	& {33065}	& {32546}	& 20	& {32154}	& {32209}	& 144	& {31963}	& {31963}	& 459	& \underline{31890}	& {31963}\\
$I_{7}$	& {32854}	& {31861}	& 17	& {31808}	& {32246}	& 116	&  \underline{31615}	& {31976}	& 316	& {31615}	& {31882}\\
$I_{8}$	& {30850}	& {30571}	& 18	& {30581}	& {30663}	& 137	& {30340}	&  \underline{30016}	& 403	& {30340}	& {30347}\\
$I_{9}$	& {34709}	& {34024}	& 17	& {34386}	& {34201}	& 162	& \underline {33911}	& {34000}	& 456	& {34197}	& {33928}\\
$I_{10}$	& {31451}	& {31201}	& 18	& {30681}	& {30891}	& 131	& {30674}	& {30811}	& 445	&  \underline{30659}	& {30674}\\
$I_{61}$	& {27158}	& {26934}	& 16	&  \underline{26644}	& {27135}	& 118	& {26644}	& {26644}	& 352	& {26644}	& {26644}\\
$I_{62}$	& {27774}	& {27619}	& 19	& {27335}	& {27342}	& 131	& {27378}	& {27293}	& 411	&  \underline{27248}	& {27248}\\
$I_{63}$	& {25308}	& {24960}	& 21	& {24518}	& {24500}	& 121	&  \underline{24312}	& {24312}	& 366	& {24312}	& {24331}\\
$I_{64}$	& {27875}	& \underline{27285}	& 21	& {27677}	& {28065}	& 143	& {27412}	& {27538}	& 462	&  \underline{27285}	& {27392}\\
$I_{65}$	& {27060}	& {26888}	& 16	&  \underline{26692}	& {27185}	& 95	& {26692}	& {26740}	& 304	& {26692}	& {26692}\\
$I_{66}$	& {27677}	& {27624}	& 19	& {27191}	& {27134}	& 133	&  \underline{27131}	& {27190}	& 392	& {27131}	& {27134}\\
$I_{67}$	& {30268}	& {30203}	& 15	& {29902}	& {29920}	& 133	& {29749}	& {30131}	& 384	&  \underline{29674}	& {29874}\\
$I_{68}$	& {28033}	& {27923}	& 24	& {27638}	& {27711}	& 163	& {27424}	& {27478}	& 430	& {27424}	&  \underline{27421}\\
$I_{69}$	& {27958}	& {27638}	& 14	& \underline {27091}	& {27091}	& 109	& {27091}	& {27148}	& 344	& {27091}	& {27425}\\
$I_{70}$	& {28483}	& {28427}	& 19	& {27930}	& {27960}	& 142	&  \underline{27789}	& {27799}	& 430	& {27789}	& {27789}\\
\hline 
\end{tabular} 
\end{center}
\begin{center}
\begin{tabular}{||c|c|c||c|c|c||c|c|c||c|c|c||}
\multicolumn{12}{l}{Tests for the set with 24 nodes visited twice}\\
\hline 
  \multirow{2}{*}& \multicolumn{2}{c||}{Results from \cite{BassettoM11}}& \multicolumn{3}{c||}{$H(3,1)$, $t_{mean}=14$s} & \multicolumn{3}{c||}{$H(5,2)$, $t_{mean}=119$s} & \multicolumn{3}{c||}{$H(6,3)$, $t_{mean}=255$s } \\ 
\cline{2-3} \cline{4-6} \cline{7-9}  \cline{10-12} 
Ins & PC & PC+m &$t_{KSH}$ & KSH & RndH &$t_{KSH}$ & KSH&RndH &$t_{KSH}$ &KSH&RndH \\ 
\hline 
$I_{41}$	& {30253}	& {30147}	& 14	& {30253}	& {30468}	& 112	& {30130}	& {30199}	& 221	&\underline{30045}	&{30161}\\
$I_{42}$	& {33008}	& {32020}	& 17	& {31787}	& {31791}	& 138	& \underline{31684}	& {31727}	& 273	& {31684}	&{31692}\\
$I_{43}$	& {31500}	& {31500}	& 15	& {31312}	& {31336}	& 128	& {31171}	&\underline {31149}	& 269	& {31149}	&{31194}\\
$I_{44}$	& {30313}	& {30170}	& 16	& {29834}	& {30068}	& 137	&\underline{29757}	& {29799}	& 269	& {29807}	&{29865}\\
$I_{45}$	& {27986}	& {27857}	& 15	& \underline{27780}	& {27839}	& 150	& {27780}	& {27843}	& 287	& {27780}	&{27780}\\
$I_{46}$	& {30073}	& \underline{30013}	& 13	& {30065}	& {30325}	& 118	& {30017}	& \underline{30013}	& 268	& {30017}	&{30092}\\
$I_{47}$	& {32106}	& {32106}	& 13	& {31694}	& {32020}	& 117	& {31694}	& {31790}	& 270	& {31694}	&\underline{31642}\\
$I_{48}$	& {31004}	& {30942}	& 15	& \underline{30470}	& {30503}	& 120	& {30470}	& {30470}	& 260	& {30470}	&{30638}\\
$I_{49}$	& {33663}	& {33185}	& 13	& {33205}	& {33411}	& 132	& {33203}	& \underline{33012}	& 307	& {33203}	&{33197}\\
$I_{50}$	& {31266}	& {31266}	& 12	& \underline{31212}	& {31236}	& 99	& {31212}	& {31212}	& 203	& {31212}	&{31212}\\
$I_{51}$	& {33722}	& {33627}	& 14	& {33358}	& {33582}	& 144	&\underline {33344}	& {33357}	& 303	& {33344}	&{33358}\\
$I_{52}$	& {32353}	& {32280}	& 14	& {32241}	& {32307}	& 97	& {32208}	& {32302}	& 226	&\underline {32199}	&{32296}\\
$I_{53}$	& {33287}	& {33200}	& 14	& {32703}	& {32871}	& 109	& \underline{32632}	& {32808}	& 224	& {32673}	&{32674}\\
$I_{54}$	& {31973}	& {31600}	& 14	& {31543}	& {31604}	& 116	& \underline{31367}	& {31368}	& 251	& {31367}	&{31525}\\
$I_{55}$	& {33837}	& \underline{33507}	& 15	& {33640}	& {33739}	& 109	& {33587}	& {33522}	& 264	& {33618}	&{33529}\\
$I_{56}$	& {29696}	& {29476}	& 13	& {29086}	& {29281}	& 106	& {29007}	& \underline{28835}	& 239	& {29003}	&{28894}\\
$I_{57}$	& {31954}	& {31640}	& 15	& {31531}	& {31585}	& 136	& {31497}	& {31615}	& 289	& \underline{31427}	&{31615}\\
$I_{58}$	& {30705}	& {30246}	& 13	& {30262}	& {30296}	& 101	& {30186}	& {30226}	& 222	&\underline{30181}	&{30559}\\
$I_{59}$	& {31549}	& {31549}	& 13	&\underline{31223}	& {31584}	& 90	& {31223}	& {31290}	& 208	& {31223}	&{31223}\\
$I_{60}$	& {32384}	& {32317}	& 15	& {32308}	& {32380}	& 116	&\underline {32136}	& {32136}	& 259	& {32136}	&{32355}\\
\hline 
\end{tabular} 
\end{center}

\end{small}

\section{Appendix: Summary of results for $50/30$ set of instances.}
 \nopagebreak[4]
{\centering
\hspace{-2.5cm}\includegraphics[scale=0.65]{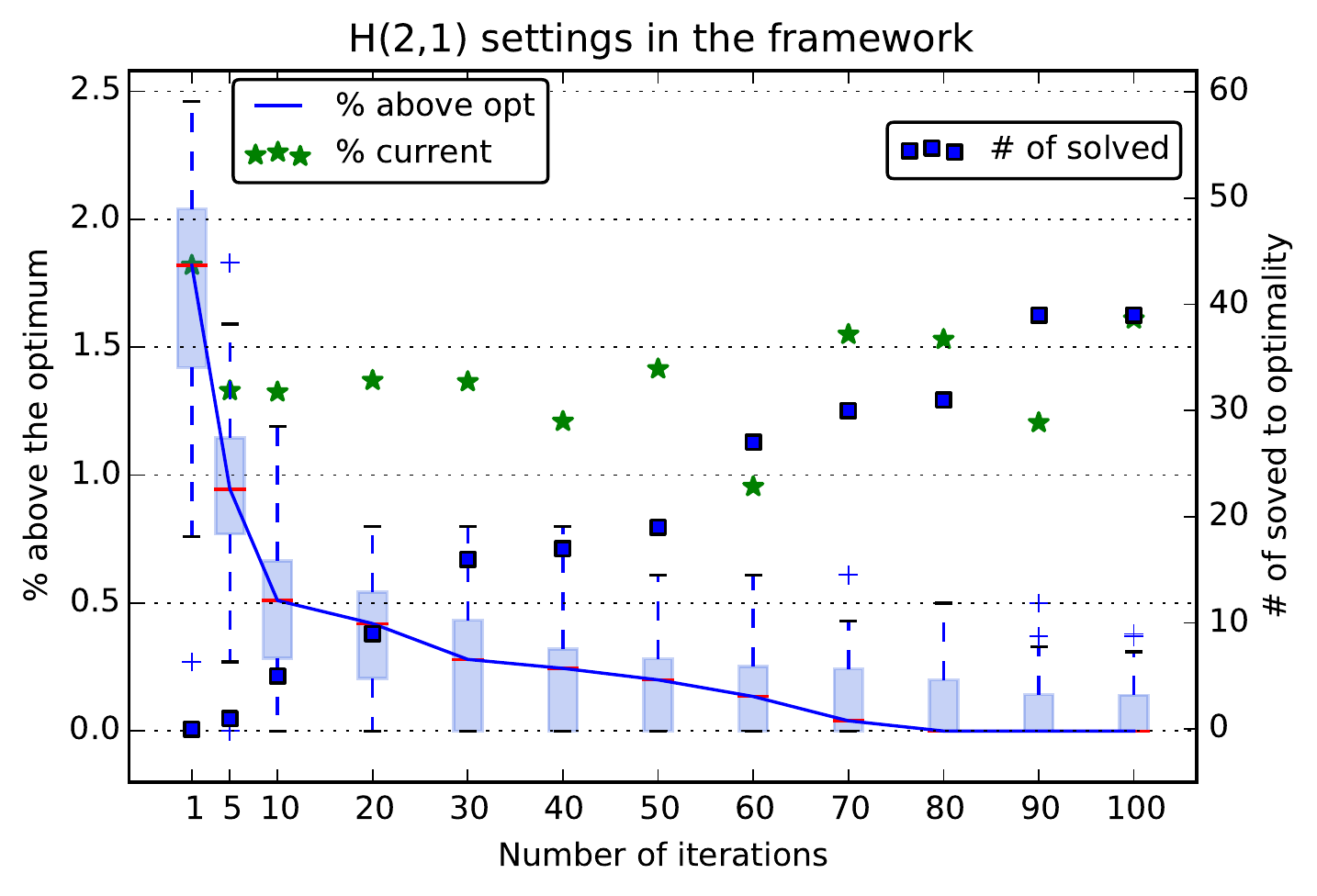}

\hspace{-2.5cm}\includegraphics[scale=0.65]{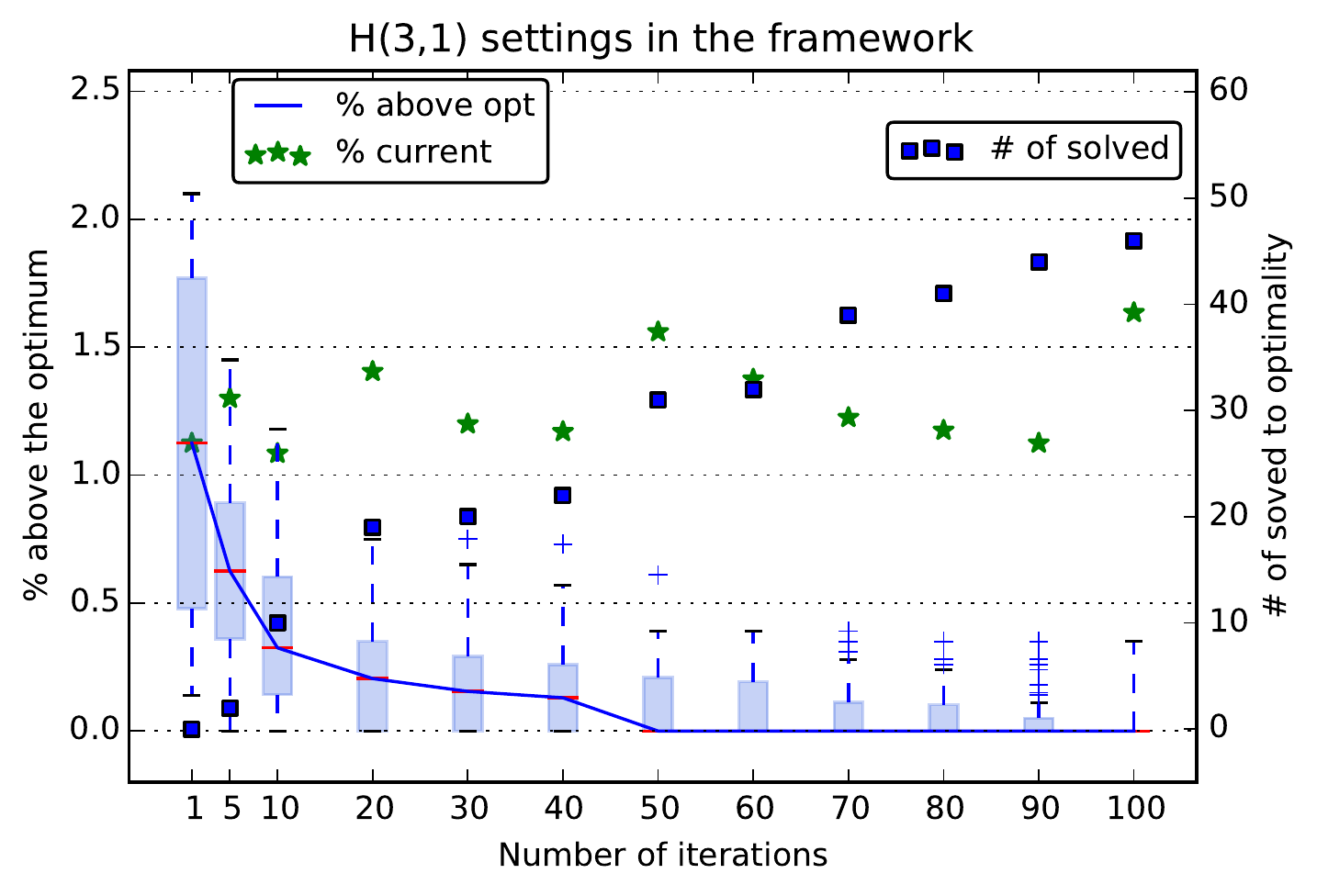}

\hspace{-2.5cm}\includegraphics[scale=0.65]{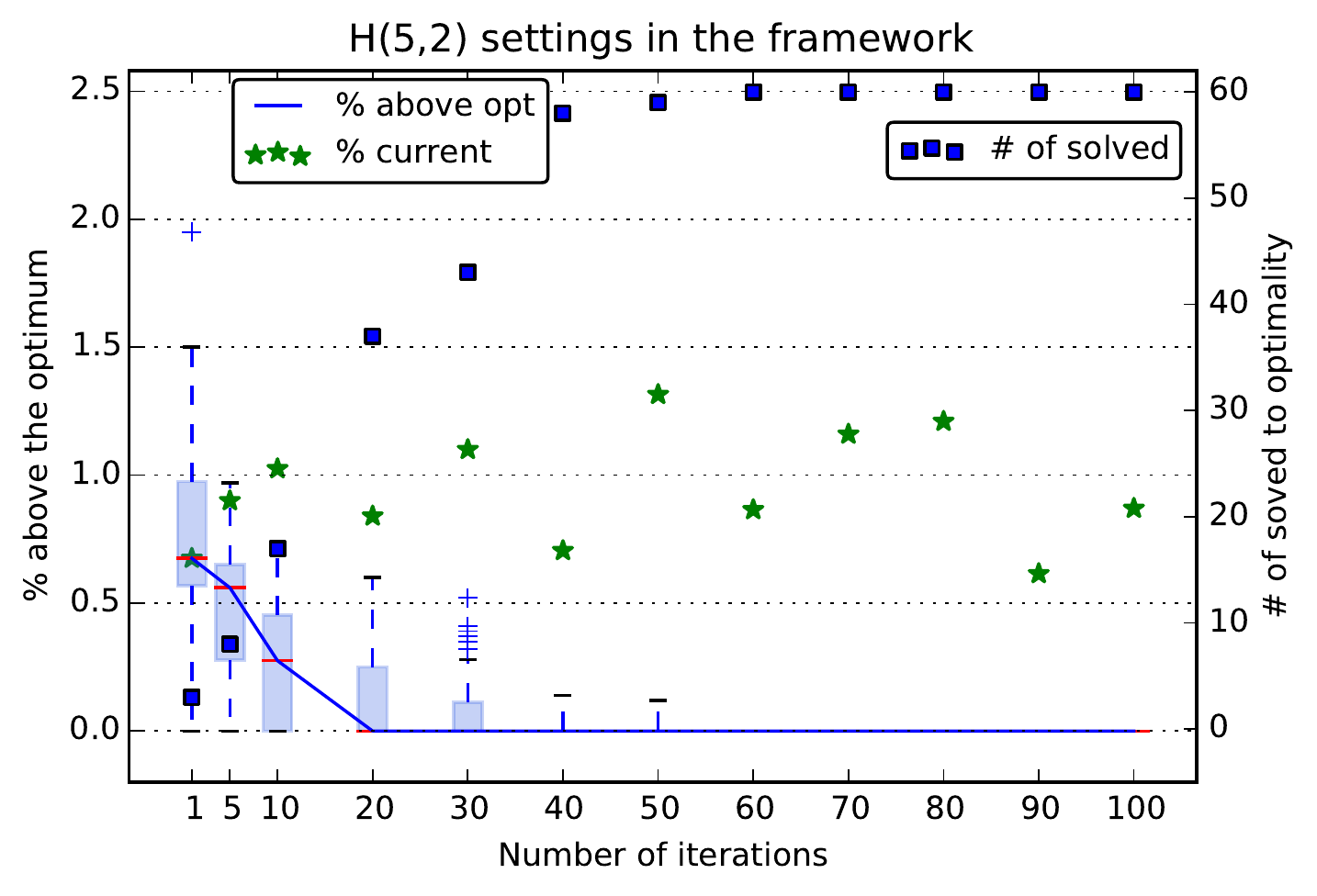}

} 

\section{Appendix: Summary of results for $100/30$ set of instances.}
 \nopagebreak[4]
{\centering
\hspace{-2.5cm}\includegraphics[scale=0.65]{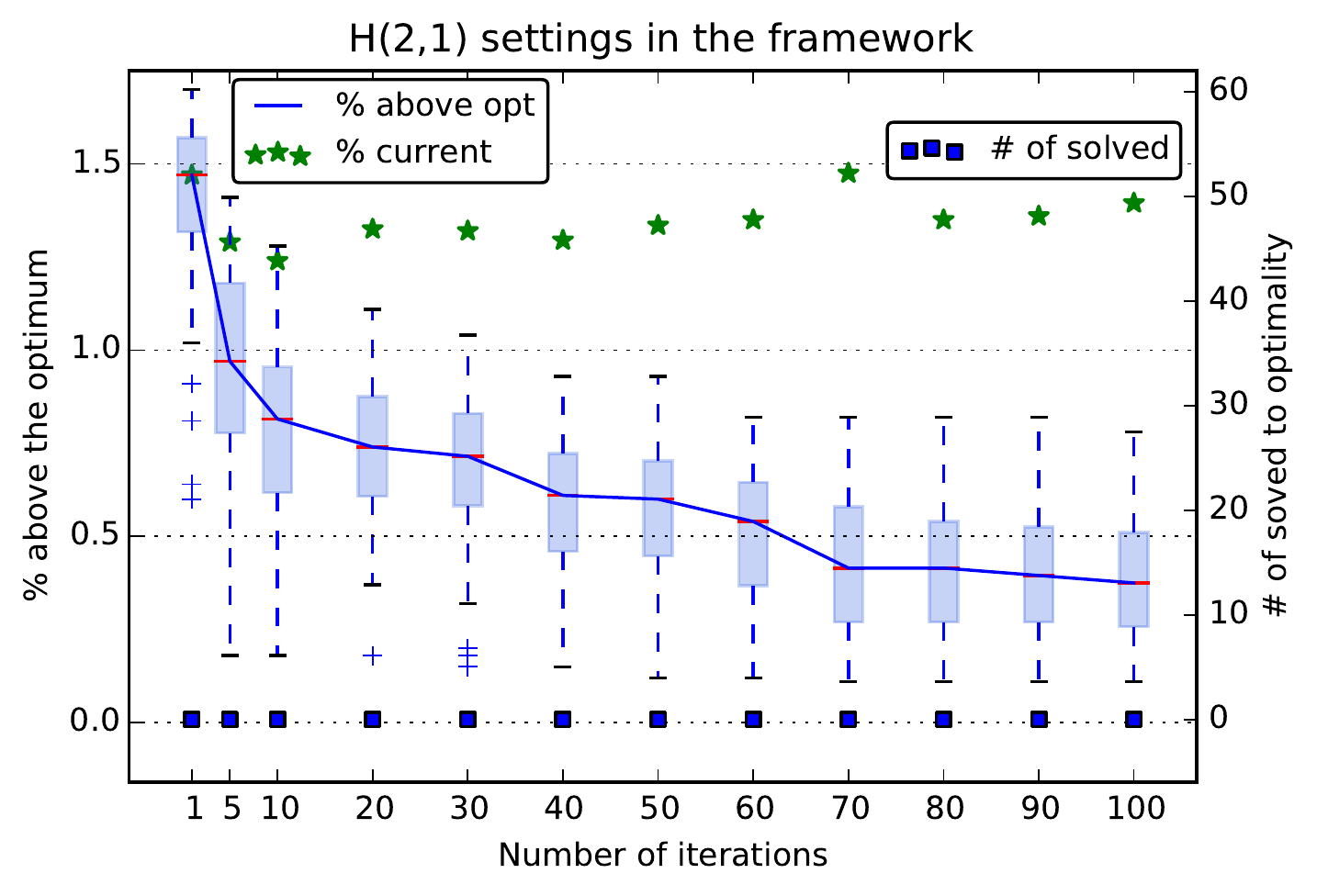}

\hspace{-2.5cm}\includegraphics[scale=0.65]{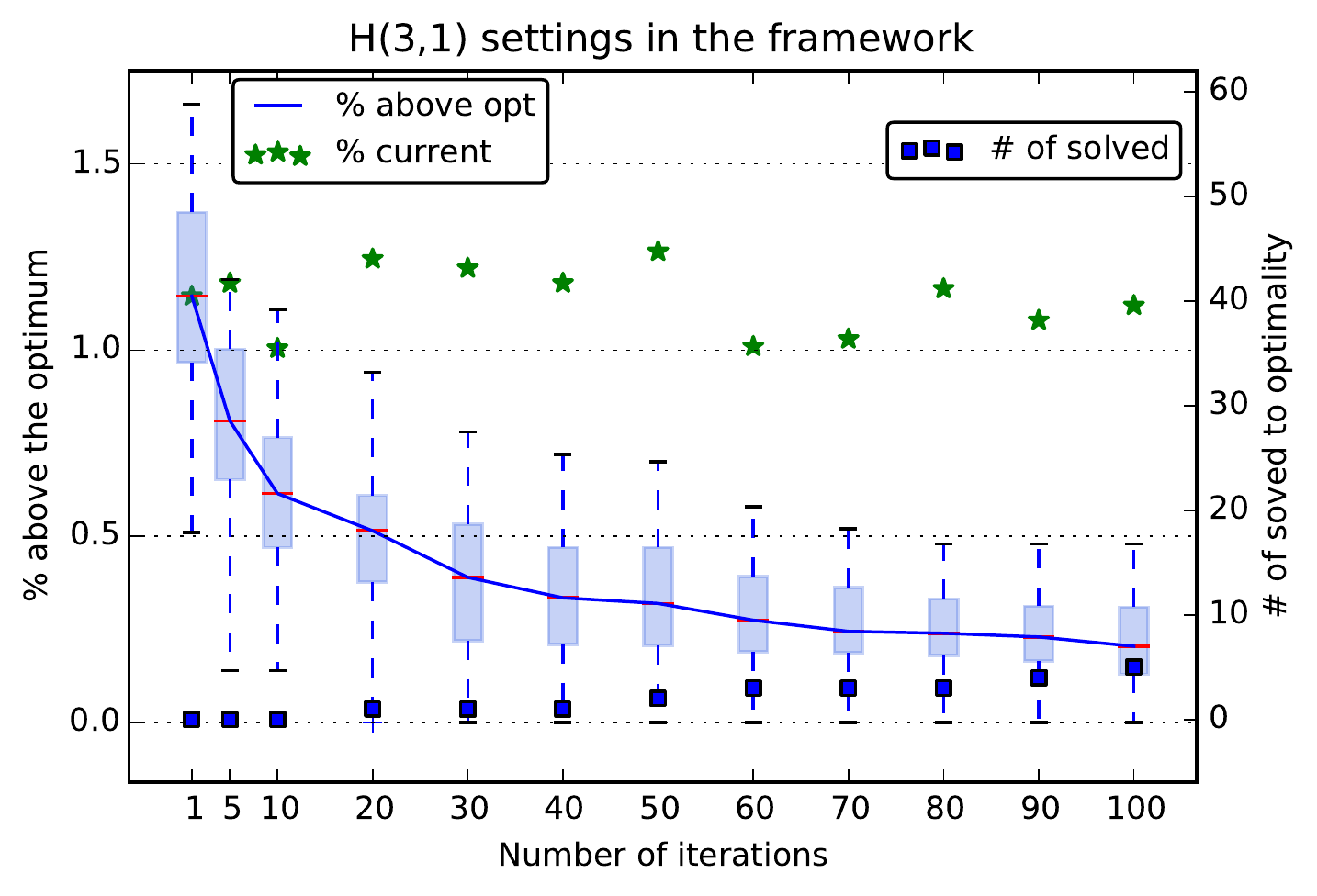}

\hspace{-2.5cm}\includegraphics[scale=0.65]{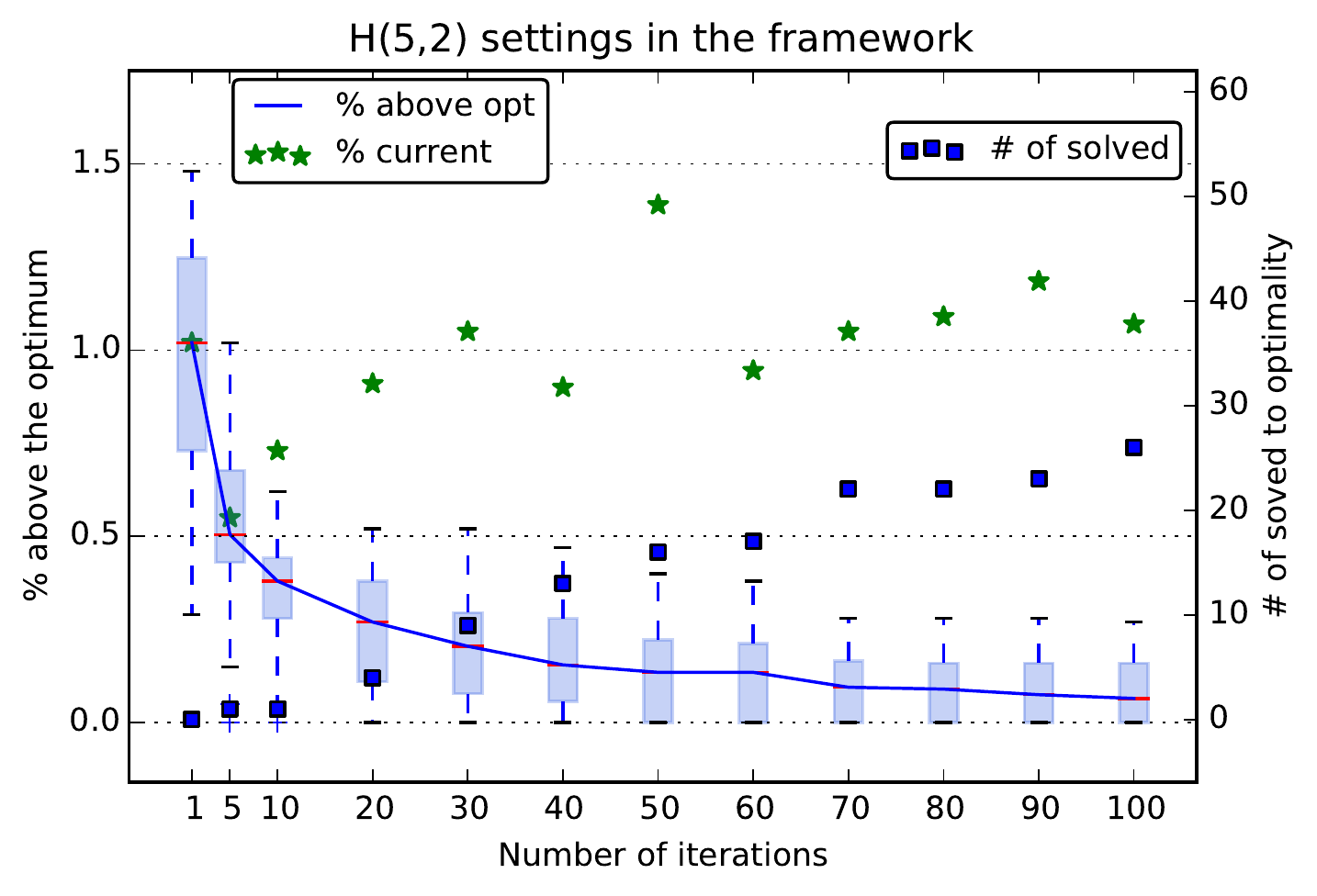}

} 


\end{document}